\documentclass[10pt]{article}

\usepackage{graphicx}
\usepackage{multicol,multirow}
\usepackage{amsmath,amssymb,amsfonts}
\usepackage{mathrsfs}
\usepackage{amsthm}
\usepackage{rotating}
\usepackage[numbers]{natbib}
\usepackage{ifpdf}
\usepackage{hyperref}
\usepackage[a4paper,left=3cm,right=2cm,top=2.5cm,bottom=2.5cm]{geometry}

\usepackage{latexsym,mathtools}
\usepackage{tikz}
\usetikzlibrary{calc}
\usepackage{soul}
\newcommand{\R}{\mathbb{R}}
\newcommand{\C}{\mathbb{C}}

\newcommand{\K}{\mathbb{K}}

\newcommand{\F}{\mathbb{F}}
\newcommand{\FR}{\mathfrak R}

\newcommand{\mathI}{\mathcal{I}}
\newcommand{\mathJ}{\mathcal{J}}

\newcommand{\wt}{\widetilde}
\newcommand{\wh}{\widehat}

\newcommand {\mat}{\left[\begin{array}}
\newcommand {\rix}{\end{array}\right]}

\newcommand{\rev}{\ensuremath{\text{rev}\,}}
\newcommand{\eqdot}{\doteq}

\newcommand{\hide}[1]{}

\def\diag{\mathop{\rm Diag}\nolimits}
\def\rank{\mathop{\rm rank}\nolimits}

\def\lcm{\mathop{\rm lcm}\nolimits}

\newtheorem{thm}{Theorem}[section]
\newtheorem{propos}[thm]{Proposition}
\newtheorem{corol}[thm]{Corollary}
\newtheorem{ex}[thm]{Example}
\newtheorem{lem}[thm]{Lemma}
\newtheorem{defi}[thm]{Definition}
\newtheorem{rem}[thm]{Remark}

\newcommand{\be}{\begin{equation}}
\newcommand{\ee}{\end{equation}}
\newcommand{\bea}{\begin{eqnarray}}
\newcommand{\eea}{\end{eqnarray}}

\makeatletter
\def\adots{\mathinner{\mkern1mu\raise\p@
\vbox{\kern7\p@\hbox{.}}\mkern2mu
\raise4\p@\hbox{.}\mkern2mu\raise7\p@\hbox{.}\mkern1mu}}
\makeatother

\title{Rosenbrock's Theorem on System Matrices over Elementary Divisor Domains}
\author{Froil\'{a}n Dopico\thanks{Universidad Carlos III de Madrid, Departamento de Matem\'aticas, Avenida de la Universidad 30, 28911 Legan\'{e}s, Spain (dopico@math.uc3m.es). Supported by grants PID2019-106362GB-I00 and RED2022-134176-T both funded by MCIN/AEI/10.13039/501100011033 and by the Madrid Government (Comunidad de Madrid-Spain) under the Multiannual Agreement with UC3M in the line of Excellence of University Professors (EPUC3M23), and in the context of the V PRICIT (Regional Programme of Research and Technological Innovation).}
\and
Vanni Noferini\thanks{Corresponding author. Aalto University, Department of Mathematics and Systems Analysis, P.O. Box 11100, FI-00076, Aalto, Finland (vanni.noferini@aalto.fi). Supported by an Academy of Finland grant (Suomen Akatemian p\"{a}\"{a}t\"{o}s 331240).}
\and
Ion Zaballa\thanks{Universidad del Pa\'is Vasco UPV/EHU, Departamento de Matem\'aticas, Apdo. Correos 644, Bilbao 48080, Spain (ion.zaballa@ehu.eus). Supported by  grant PID2021-124827NB-I00 funded by MCIN/AEI/10.13039/501100011033 and by “ERDF A way of making Europe” by the “European Union”.}}

\begin{document}
\maketitle
\abstract{Rosenbrock's theorem on polynomial system matrices is a classical result in linear systems theory that relates the Smith-McMillan form of a rational matrix $G$ with the Smith forms of an irreducible polynomial system matrix $P$ giving rise to $G$ and of a submatrix of $P$. This theorem has been essential in the development of algorithms for computing the poles and zeros of a rational matrix via linearizations and generalized eigenvalue algorithms. In this paper, we extend Rosenbrock's theorem to  system matrices $P$ with entries in an arbitrary elementary divisor domain $\FR$ and matrices $G$ with entries in the field of fractions of $\FR$. These are the most general rings where the involved Smith-McMillan and Smith forms both exist and, so, where the problem makes sense. Moreover, we analyze in detail what happens when the system matrix is not irreducible. Finally, we explore how Rosenbrock's theorem can be extended when the system matrix {$P$} itself has entries in the field of fractions of the elementary divisor domain.}



\section{Introduction} In the classical monograph \cite{Rosen70}, Rosenbrock studies very general multivariable
input-output systems defined by linear algebraic  and differential equations. After taking the Laplace transform with zero initial conditions, the coefficients of the obtained system \cite[Chapter 2]{Rosen70} can be arranged  in a single polynomial matrix, the  so-called \textit{polynomial system matrix}: $P(s) = \begin{bsmallmatrix}
A(s) & B(s) \\ C(s) & D(s) \end{bsmallmatrix}$  where $A(s)$ is always nonsingular. Then Rosenbrock observed that the matrix connecting the Laplace transforms of the inputs and the outputs, known as the \emph{transfer function} matrix of the system, is the Schur complement of $A(s)$ in $P(s)$, that is, $G(s) = D(s) - C(s) A(s)^{-1} B(s)$. This is a rational matrix, that is, a matrix whose entries are univariate rational functions with coefficients in a field. It is also shown in \cite[Section 5, Chapter 3]{Rosen70} that any rational matrix $G(s)$ can be represented, or realized, by a quadruple of polynomial matrices $A(s)$, $B(s)$, $C(s)$ and $D(s)$ as $G(s) = D(s) - C(s) A(s)^{-1} B(s)$, where $A(s)$ is nonsingular. The notion of minimality is of fundamental importance in many results of \cite{Rosen70}.
A polynomial system matrix is said to be \emph{minimal}, or \emph{of least order}, or \emph{irreducible}, if $A(s)$ and $B(s)$ are left coprime and $A(s)$ and $C(s)$ are right coprime\footnote{For more details on left and right coprimeness, and other concepts mentioned in this introduction, see Section \ref{sec.background}.}. Among the many results proved by Rosenbrock on irreducible polynomial system matrices, this paper focuses on \cite[Theorem 4.1, Chapter 3]{Rosen70}, which relates the pole and zero structures of the rational matrix $G(s)$ with the zero structures of $A(s)$ and $P(s)$, respectively. In other words, \cite[Theorem 4.1, Chapter 3]{Rosen70}
relates the Smith-McMillan form {\ \cite[p. 109]{Rosen70}} of the rational matrix $G(s)$ with the Smith forms {\ \cite[p. 8]{Rosen70}} of $A(s)$ and $P(s)$ (see Theorem \ref{thm.1stproof}). In this paper, we refer to this result as \emph{Rosenbrock's theorem}; it was originally proved by Rosenbrock for system matrices over the ring of univariate polynomials with coefficients in the real or complex fields \cite[Theorem 4.1, Chapter 3]{Rosen70}.

Rosenbrock's theorem is particularly important in the numerical computation of the poles and the zeros of a rational matrix. This is due to the combination of two facts: first, for every rational matrix $G(s)$ there exist {\em linear} irreducible polynomial system matrices whose Schur complements are $G(s)$ and, second, there exist efficient and reliable algorithms for computing the eigenvalues (the zeros) of linear polynomial matrices \cite{moler-stewart,vd79}. This approach was thoroughly explored in the classical paper \cite{vd81}. Moreover, the numerical computation of the zeros of rational matrices has received considerable attention in recent years since rational matrices have appeared in new applications \cite{mehrmann-voss,subai11} and also in the approximation of more general nonlinear eigenvalue problems \cite{guttel-tisseur}. The rational matrices arising in these new applications and approximations are often represented in different ways than those arising in linear systems theory and this has motivated an intense research in the construction of linear irreducible polynomial system matrices related to such representations. See \cite{alambehera,AmDoMaZa18}, among many other references on this topic. Rosenbrock's theorem has played a key role in all these recent developments.

It is known that many of Rosenbrock's results on rational matrices and polynomial system matrices, and in particular Rosenbrock's theorem, are valid for system matrices over an arbitrary principal ideal domain and their transfer matrices, having entries in the field of fractions. This was first proved by Coppel in \cite{Coppel74}, just a few years after the publication of \cite{Rosen70}. Coppel followed a considerably more elementary approach than the one in \cite{Rosen70}. Unfortunately, \cite{Coppel74} has not received all the attention that it deserves. For instance, Rosenbrock's theorem implies that the computation of the Smith form of a large integer matrix can, under somewhat mild assumptions, be reduced to the computation of the Smith-McMillan form of a (potentially much smaller) matrix of rational numbers. To the best of our knowledge, this fact has not been explored in the development of algorithms for integer matrices.

The main contribution of this paper is to extend Rosenbrock's theorem  to matrices with entries in the field of fractions of any elementary divisor domain and to their irreducible system matrices.  Elementary divisor domains were introduced by Kaplansky \cite{kap49}. They are precisely those integral domains where every matrix  is unimodularly equivalent to one in Smith normal form. As a consequence, every matrix with entries in the field of fractions of an elementary divisor domain is unimodularly equivalent to one in Smith-McMillan form. Therefore, it is easy to see that elementary divisor domains are the most general rings where the statement of Rosenbrock's theorem makes sense and may be true. Hence, in this paper we show that Rosenbrock's theorem is true under the mildest conceivable assumptions on the base ring. A very important example of an elementary divisor domain\footnote{A good general reference for elementary divisor domains and other related integral domains is \cite[Chapter 1]{Fried16}. See also Appendix \ref{append} in this paper.} that is not a principal ideal domain (PID) is the ring of complex-valued functions that are holomorphic on an open\footnote{Openness is crucial: in contrast, the ring of functions that are holomorphic on a \emph{compact} connected set is always a PID.} connected subset of $\mathbb{C}$ \cite{hel40,hel43}. We believe that Rosenbrock's theorem for these holomorphic functions might have applications in the numerical solution of some nonlinear eigenvalue problems. Other interesting examples of elementary
divisor domains that are not principal ideal domains include the ring of algebraic integers \cite[Remark {17}]{N23},
or  the ring $\{p\in\R(s)[z,z^{-1}]: p(s,e^{-s}) \text{ is holomorphic in } \C\}$ \cite{CaQue15,Gl-Lu97}.

A second contribution of this work is a detailed study of the information that can be extracted from a \emph{reducible} (that is, not minimal) system matrix about the Smith-McMillan form of its Schur complement. This study has practical implications in situations where linear polynomial system matrices which are not guaranteed to be irreducible are used for the numerical computation of the zeros of a rational matrix \cite{nleigs}.

An additional contribution of this work is a  generalization of Rosenbrock's theorem to the case where the system matrix itself has entries in the field of fractions of the elementary divisor domain. In the case of rational functions with coefficients in a field, this means that in $G(s) = D(s) - C(s) A(s)^{-1} B(s)$ the matrices $A(s), B(s), C(s)$ and $D(s)$ are rational matrices. As far as we know, this extension has never been considered even for rational functions and, apart from its fundamental interest, it may have applications for computing the poles and zeros at infinity of univariate rational matrices.

Finally, we emphasize that, in addition to its weaker assumptions on the base ring, the proof presented in this paper of Rosenbrock's theorem is more direct than other proofs available in the literature. For instance, the original one presented by Rosenbrock in \cite{Rosen70} relies on nontrivial properties of the strict system equivalence relation of polynomial system matrices. Moreover, this implies that the proof in \cite{Rosen70} is not valid  when the dimension of $A(s)$ in $P(s) = \begin{bsmallmatrix} A(s) & B(s) \\ C(s) & D(s) \end{bsmallmatrix}$ is smaller than the degree of $\det A(s)$. The classical reference \cite{Kail80} also includes in its Section 8.3 a proof of Rosenbrock's theorem over the ring of univariate polynomials with coefficients in the complex field (though the result itself is not clearly stated as a theorem), which is again based on different concepts of equivalence of polynomial system matrices. The proof in \cite{Coppel74} holds over a PID and relies on the properties of the so-called determinantal denominators of matrices with entries in the field of fractions of a principal ideal domain. Some of these properties are not trivial and to extend them to arbitrary elementary divisor domains appears challenging.

This paper is organized as follows. Section \ref{sec.background} includes preliminaries and background material. Section \ref{sec.1stproof} presents the proof of Rosenbrock's theorem over elementary divisor domains. Section \ref{sec.beyond} explores what happens when the system matrix is not irreducible. Three interesting corollaries of Rosenbrock's theorem are presented in Section \ref{sec.consequences} and, in Section \ref{sec.RTFF},  Rosenbrock's theorem is generalized to the case
when the involved matrices have their entries in the field of fractions of an elementary divisor domain.  Conclusions are discussed in Section
\ref{sec.conclusions} and  two appendices present proofs valid under weaker assumptions of certain technical results that are well known when the underlying ring is a principal ideal domain.

\section{Background material} \label{sec.background}
Throughout this paper $\FR$ denotes an elementary divisor domain (EDD) \cite[Section 1.5]{Fried16}\footnote{To avoid having to pedantically exclude trivial special cases in our statements, we will also assume throughout that $\FR$ is not a field and is not the trivial ring $\{ 0 \}$. In particular, there is $r \in \FR$ such that $0 \neq r \not \in {\FR^\times}$.},
$\F$ denotes the field of fractions of $\FR$ and ${\FR^\times}$ denotes the group of units of $\FR$.
It is useful to bear in mind that every EDD is a B\'{e}zout domain (BD) and a greatest common divisor domain (GCDD; see Appendix \ref{append}).
In particular, in a GCDD $\FR$ every pair of elements $a,b\in\FR$ have a greatest
common divisor, denoted by $\gcd(a,b)$, and a least common multiple, denoted by $\lcm (a,b)$. Moreover, they satisfy
$\lcm (a,b)=\frac{a b}{\gcd(a,b)}$. Both $\gcd(a,b)$ and $\lcm (a,b)$ are defined up to products by units of $\FR$. In addition, in a BD $\gcd(a,b)$ is always an $\FR$-linear combination of $a$ and $b$, i.e., there exist $x,y \in \FR$ such that $\gcd(a,b)=ax+by$.
EDDs are characterized by the property that, for all triples $a,b,c\in \FR$, there exist $x,y,z,w \in \FR$ such that $\gcd(a,b,c)=zxa+zyb+wyc$ \cite[{Definition 1.5.1}]{Fried16}. For some auxiliary results, it suffices to assume that $\FR$ is a BD or a GCDD: when this is the case, we will emphasize it explicitly by only making the weaker assumption in the statement. On the other hand, not always an EDD is a unique factorization domain (UFD; see Appendix \ref{append}); more precisely, an EDD is a UFD if and only if it is a principal ideal domain (PID). As a consequence, any argument that uses factorizations of elements into products of irreducible or prime elements is not valid when working with EDDs that are not PIDs. This complicates the proofs of some results that are well known in the standard literature on polynomial and rational matrices since they are often, either explicitly or implicitly, based on arguments that assume the existence of factorizations into prime factors{; see, for example,} \cite{Coppel74,Kail80,Rosen70}.

The sets of $p\times m$ matrices with entries in $\FR$ and $\F$ are denoted by $\FR^{p\times m}$ and $\F^{p\times m}$, respectively. A \emph{unimodular} matrix $A\in \FR^{n\times n}$ is a square matrix which is invertible over $\FR$, that is, there exists a matrix $A^{-1}\in \FR^{n\times n}$ such that $AA^{-1} = A^{-1} A =I_n$, where $I_n$ is the identity matrix of size $n\times n$. Clearly, $A\in \FR^{n\times n}$ is unimodular if and only if its determinant is a unit, i.e., $\det A \in {\FR^\times}$ \cite[Lemma 1.7.6]{Fried16}. A square $n \times n$ matrix $A$ with entries in $\F$ (or in $\FR$) such that $\det A \ne 0$ is said to be nonsingular. Equivalently, a square matrix is nonsingular if it is invertible over $\F$; note that being unimodular implies being nonsingular, but the converse is not generally true. The rank of a $p\times m$ matrix $A$ with entries in $\F$ (or in $\FR$) is the maximal size of a nonzero minor of $A$ and is denoted by $\rank A$.

Square diagonal matrices arise very often in this paper. Therefore, we use a shorthand notation for them: $ \diag \left(\alpha_1,\ldots, \alpha_r\right)$ denotes the $r\times r$ diagonal matrix whose diagonal entries are $\alpha_1,\ldots, \alpha_r$. Block diagonal matrices, or direct sums of matrices, also arise very often and they are denoted in two different ways as
\[
\begin{bmatrix}
  A & 0 \\
  0 & B
\end{bmatrix} =: A\oplus B \in \F^{(m+p) \times (n+q)},
\]
where $A \in \F^{m\times n}$ and $B \in \F^{p\times q}$.

The zero matrix of size $p\times q$ is denoted by $0_{p\times q}$ or simply by $0$, when the size is clear from the context or is not relevant. To simplify the statements of some results the degenerate cases $0_{0\times 0}$, $0_{0\times q}$ and $0_{p\times 0}$ are also used. $0_{0\times 0}$ stands for the empty matrix, $A \oplus 0_{0\times q}$ indicates that $q$ zero columns are appended to the right of $A$, while $A \oplus 0_{p\times 0}$ indicates that $p$ zero rows are appended to the bottom of $A$.

Two matrices $A,B \in \F^{p\times m}$ are called unimodularly equivalent, or (for brevity) simply \emph{equivalent} in this paper, if there exist unimodular matrices $U\in \FR^{p\times p}$ and $V\in \FR^{m\times m}$ such that $UAV = B$. We emphasize that $A$ and $B$ are allowed to possibly have entries in $\F$, although the definition also applies to matrices over $\FR$, while $U$ and $V$ are required to have entries in $\FR$.

Every matrix with entries in an EDD has a Smith form (see, for example, \cite[Theorem 1.14.1]{Fried16} or \cite[Chapter 2]{Shchedrykbook}).
That is, for each $A \in \FR^{p\times m}$ with $\rank A = r$, there exist unimodular matrices $U\in \FR^{p\times p}$ and $V\in \FR^{m\times m}$ such that
\begin{equation} \label{eq.smith}
U A V = \begin{bmatrix} \diag \left(\alpha_1,\ldots, \alpha_r\right)&0\\0&0
\end{bmatrix} =: S_A,
\end{equation}
where $\alpha_i \ne 0$ for $i =1,\ldots r$ and $\alpha_1\mid \alpha_2 \mid \cdots\mid
\alpha_{r}$.
This was first proved by Smith \cite{Smith} for $\FR=\mathbb{Z}$ and extended
to different classes of commutative and non-commutative rings by several authors (see
the preface of \cite{Shchedrykbook} and the remarks at the end of Chapter 1 of \cite{Fried16} for a brief historical account). A proof for any EDD was given by Kaplansky \cite{kap49}.

The matrix $S_A \in \FR^{p\times m}$ is a Smith form of $A$ and $\alpha_1, \ldots , \alpha_r$ are called \emph{invariant factors of $A$}. The invariant factors of $A$ are uniquely determined by $A$ up to a product by units of $\FR$. Equivalently, the Smith form of $A$ is uniquely determined by $A$ up to products by unimodular diagonal matrices. Thus, we have that $A$, $B \in \FR^{p\times m}$ are equivalent if and only if $A$ and $B$  share a common Smith form. We will say that an invariant factor of $A$ is nontrivial when it is not a unit of $\FR$.

Many of the quantities of interest in this paper are only defined up to multiplication by a unit of $\FR$, including for example invariant factors or greatest common divisors of sets of elements. Thus the equations or equalities where such quantities are involved should rigorously include a unit of $\FR$. To avoid a frequent repetition of this fact, we will simply use the symbol $\doteq$ to indicate that a scalar equality holds ``up to a product'' by an element of ${\FR^\times}$. In the same spirit, we will speak about ``the" greatest common divisor of some elements or ``the" invariant factors of a matrix, even though strictly speaking these are defined as equivalence classes of associate elements of $\FR$. Analogously, we will speak about ``the" Smith form of a given matrix $A \in \FR^{p \times m}$ to mean any of the matrices of the form \eqref{eq.smith} that are equivalent to $A$ and that differ one from another by (left or right) multiplication times a diagonal unimodular matrix. Moreover, we will also use $\doteq$ for matrix equalities that hold ``up to a (left or right) multiplication by a unimodular diagonal matrix". For example, if the $r$ invariant factors of $A \in \FR^{p\times m}$ in \eqref{eq.smith} are all trivial, we may write that ``the" Smith form of $A$ is $S_A \doteq I_r \oplus 0$.

The invariant factors of $A\in\FR^{p\times m}$ can be determined
in terms of the elements of $A$ as follows \cite[{Lemma 1.12.3}]{Fried16}: for $k=1,\ldots, \min\{p,m\}$, let $\delta_k(A)$ be the greatest common divisor of the minors of order $k$ of $A$. Then, $\delta_k(A)\mid\delta_{k+1}(A)$, $k=1,\ldots, \min\{p,m\}-1$ and
\be\label{eq.ifdetdiv}
\delta_k(A) \doteq
\alpha_1\alpha_2\cdots \alpha_k,\quad  k=1,\ldots, r.
\ee
Thus, %up to a product by a unit of $\FR$,
$\alpha_k \doteq \frac{\delta_k(A)}{\delta_{k-1}(A)}$, $k=1,\ldots, r$ where
$\delta_0(A)=1$. The elements $\delta_1(A)\mid\delta_2(A)\mid\cdots\mid\delta_{\min\{p,m\}}(A)$ are called the \textit{determinantal divisors}
of $A$.

Every matrix with entries in the field of fractions of an EDD has a so-called Smith-McMillan form; this fact was first established by McMillan \cite{McM} for the case $\FR=\C[z]$ and $\F=\C(z)$. Namely, for every $G \in \F^{p\times m}$ with $\rank G = r$, there exist unimodular matrices $U\in \FR^{p\times p}$ and $V\in \FR^{m\times m}$ such that
\begin{equation} \label{eq.smithmcmillan}
U G V = \begin{bmatrix} \diag \left(\frac{\varepsilon_1}{\psi_1},\ldots, \frac{\varepsilon_r}{\psi_r}\right)&0\\0&0
\end{bmatrix} =: S_G,
\end{equation}
where $0\ne \varepsilon_i \in \FR$ and $0 \ne \psi_i \in \FR$ for $i=1, \ldots, r$, $\varepsilon_1\mid\varepsilon_2\mid\cdots\mid \varepsilon_r$,
$\psi_r\mid\psi_{r-1}\mid\cdots\mid \psi_1$ and $\gcd(\varepsilon_i,\psi_i){ \doteq 1}$, $i=1,\ldots, r$.
These elements of $\FR$ are uniquely determined by $G$ up to a product by units of $\FR$ and
$\frac{\varepsilon_1}{\psi_1},\ldots, \frac{\varepsilon_r}{\psi_r}$ are called the \textit{invariant fractions} of $G$.
The matrix $S_G \in \F^{p\times m}$ is a Smith-McMillan form of $G$ and is again uniquely determined by $G$ up to products by unimodular diagonal matrices. Thus, we have that $G$, $H \in \F^{p\times m}$ are equivalent if and only if $G$ and $H$ share a common Smith-McMillan form.

The existence of the Smith-McMillan form of every $G \in \F^{p\times m}$ can be obtained as a corollary of the existence of the Smith form for matrices over $\FR$ as follows: let $\varphi$ be a least common denominator of the entries of $G$. Then $\varphi G\in\FR^{p\times m}$ and thus it has a Smith form
$
\diag \left(\alpha_1,\ldots, \alpha_r\right) \oplus 0_{(p-r) \times (m-r)}$.
Then cancelling the greatest common divisor of $\varphi$ and $\alpha_i$, for $i=1, \ldots, r$, leads to the Smith-McMillan form $S_G \in\F^{p\times m}$ in \eqref{eq.smithmcmillan}. Proving this fact is not difficult when $\FR$ is a UFD; this is the case when $\FR$ is a PID \cite{Coppel74}.  To prove it more generally when $\FR$ is an EDD, but not a PID, we can use Lemma \ref{lem_divnumden} in Appendix \ref{append}.
Applying that lemma to  $\alpha_1\mid\alpha_2\mid\cdots\mid \alpha_r$ and $b_1=b_2=\ldots= b_r=\varphi$ we
obtain $\epsilon_i\doteq \frac{\alpha_i}{\gcd(\alpha_i,\varphi)}\;\big|\;\frac{\alpha_{i+1}}{\gcd(\alpha_{i+1},\varphi)}\doteq\epsilon_{i+1}$,
and $\psi_{i+1}\doteq\frac{\varphi}{\gcd(\alpha_{i+1},\varphi)}\;\big|\;\frac{\varphi}{\gcd(\alpha_{i},\varphi)}\doteq\psi_{i}$ for $i=1,\ldots, r-1$.

Let $G_1\in\FR^{p\times m}$ and $G_2\in\FR^{q\times m}$. A matrix $D \in\FR^{m\times m}$ is a common right divisor of $G_1$ and $G_2$ if there exist matrices $\widetilde{G_1}\in\FR^{p\times m}$ and $\widetilde{G_2}\in\FR^{q\times m}$ such that
\begin{equation} \label{eq.defcommondivisor}
G_1 = \widetilde{G_1} \, D \quad \mbox{and} \quad G_2 = \widetilde{G_2} \, D.
\end{equation}
The matrices $G_1$ and $G_2$ are said to be right coprime if every common right divisor is unimodular. Two matrices $H_1\in\FR^{n\times p}$ and $H_2\in\FR^{n\times q}$ are left coprime if $H_1^T$ and $H_2^T$ are right coprime.

It is known that if $p+q < m$, then $G_1\in\FR^{p\times m}$ and $G_2\in\FR^{q\times m}$ cannot be right coprime because they always have non-unimodular right divisors \cite[Remark 3.1]{AmMaZa16}. In the case $p+q \geq m$ it is not difficult to prove the characterizations of coprime matrices over EDDs in Proposition \ref{teo_copr_an}. Note that such characterizations are well known when $\FR$ is the ring of univariate polynomials with coefficients in any field; see, for instance, \cite{Kail80, NV23,Rosen70}. For completeness, a proof of Proposition \ref{teo_copr_an} valid over EDDs is included in Appendix \ref{appendcoprimeness}.

\begin{propos}\label{teo_copr_an}
Let $G_1\in\FR^{p\times m}$ and $G_2\in\FR^{q\times m}$, $p+q\geq m$.
The following are equivalent:
\begin{description}
\item[\rm i)] $G_1$ and $G_2$ are right coprime in $\FR$.
\item[\rm ii)] The Smith form over $\FR$ of $\begin{bsmallmatrix}G_1\\G_2\end{bsmallmatrix}$
is $\begin{bsmallmatrix}I_m\\0\end{bsmallmatrix}$.
\item[\rm iii)] There exists a unimodular matrix $U\in\FR^{(p+q)\times (p+q)}$ such that
$U\begin{bsmallmatrix}G_1\\G_2\end{bsmallmatrix}=\begin{bsmallmatrix}I_m\\0\end{bsmallmatrix}$.
\item[\rm iv)] There exist matrices $C\in\FR^{p\times (p+q-m)}$, $D\in\FR^{q\times (p+q-m)}$
 such that $\begin{bsmallmatrix}G_1&C\\G_2&D\end{bsmallmatrix}$ is unimodular.
\item[\rm v)] There exist matrices $X\in\FR^{m\times p}$, $Y\in\FR^{m\times q}$ such that $XG_1+YG_2=I_m$.
\end{description}
\end{propos}

An analogous result can be given for left coprime matrices with elements in $\FR$ (simply transpose everything in Proposition \ref{teo_copr_an}). The characterizations in Proposition \ref{teo_copr_an} and their counterparts for left coprime matrices will be often used without being explicitly referred to.

\section{Rosenbrock's theorem over EDDs} \label{sec.1stproof}
This section uses some results developed by Coppel in  \cite{Coppel74} for PIDs. However, we weaken the assumption about the base ring and consider more general EDDs. Moreover, in stark contrast with \cite{Coppel74}, we do not use determinantal denominators, which play a fundamental role in Coppel's proof of the ``denominator part'' of Rosenbrock's theorem, as well as for several side results.
Rosenbrock's theorem over EDDs is proved in Theorem \ref{thm.1stproof}. Its proof relies on {some auxiliary results:} {Proposition} \ref{thm.1stpraux1}, {Theorem} \ref{thm.1stpraux2} and {Theorem} \ref{thm.1stpraux3}. {Proposition} \ref{thm.1stpraux1} can be found (for a PID)  in \cite[p. 101]{Coppel74}, though not formally stated as a theorem. For the case $\FR=\C[x]$ and $\F=\C(x)$ (univariate polynomial  and rational matrices with complex coefficients) it also appears in \cite[Theorem 6.5-4]{Kail80}, with a different proof. Our proof is inspired by that in \cite{Coppel74} but holds in broader generality when $\FR$ is any EDD.

\begin{propos} \label{thm.1stpraux1} Let $A \in \FR^{n\times n}$, with $\det A \ne 0$, $C \in \FR^{p\times n}$ and $R = C A^{-1}\in \F^{p\times n}$. Assume that $A$ and $C$ are right coprime. Then:
\begin{description}
  \item[\rm i)] If $Y \in \FR^{n\times q}$ satisfies that $R Y \in \FR^{p\times q}$, then there exists a matrix $X \in \FR^{n\times q}$ such that $Y = A X$.
  \item[\rm ii)] If $\widetilde{A} \in \FR^{n\times n}$ and $\widetilde{C} \in \FR^{p\times n}$ satisfy the same properties of $A$ and $C$, i.e., (1) $\det \widetilde{A} \neq 0$ (2) $\widetilde{A}$ and $\widetilde{C}$ are right coprime and (3) $R = \widetilde{C} \widetilde{A}^{-1}$, then there exists a unimodular matrix $H \in \FR^{n\times n}$ such that $\widetilde{A} = A H$ and $\widetilde{C} = C H$.
  \end{description}
\end{propos}
\smallskip \noindent
\textbf{Proof}:
\begin{description}
\item[\rm i)] Since $A$ and $C$ are right coprime, there exist $M\in \FR^{n\times n}$ and $N \in \FR^{n\times p}$ such that $MA + NC = I_n$. Thus, $A^{-1} = M + NC A^{-1}= M + N R$. This implies $A^{-1} Y = MY + N (RY) \in \FR^{n\times q}$. Picking $X := (M+NR)Y$ concludes the proof.
\item[\rm ii)] Since $R \widetilde{A} = \widetilde{C} \in \FR^{p\times n}$, item i) implies that there exists $H \in \FR^{n\times n}$ such that $\widetilde{A} = A H$. Therefore $\widetilde{C} = R \widetilde{A} = (R A) H = CH$. But $\widetilde{A}$ and $\widetilde{C}$ are right coprime, and hence $H$ is unimodular. \qed \smallskip
\end{description}

Theorem \ref{thm.1stpraux2} below can be seen as Rosenbrock's theorem for a particular class of ``realizations'' of matrices with entries in the field of fractions of an EDD. The second part of the proof of Theorem \ref{thm.1stpraux2} relies on an argument that can be found at the bottom of \cite[p. 105]{Coppel74}; the first part follows from {Proposition} \ref{thm.1stpraux1}. For univariate polynomial and rational matrices with complex coefficients{,} related results can be found in \cite[p. 446]{Kail80}.

\begin{thm} \label{thm.1stpraux2} Let $A \in \FR^{n\times n}$, with $\det A \ne 0$, $C \in \FR^{p\times n}$ and $R = C A^{-1}\in \F^{p\times n}$. Assume that $A$ and $C$ are right coprime. If the Smith-McMillan form of $R$ is
$$
S_R \doteq \diag \left(\frac{\varepsilon_1}{\psi_1},\ldots, \frac{\varepsilon_r}{\psi_r}\right) \oplus 0_{(p-r) \times (n-r)} \in \F^{p\times n},
$$
then
\begin{description}
  \item[\rm i)] the Smith form of $C$ is $S_C \doteq \diag \left(\varepsilon_1,\ldots, \varepsilon_r \right) \oplus 0_{(p-r) \times (n-r)} \in \FR^{p\times n},$ and
  \item[\rm ii)] the Smith form of $A$ is $S_A \doteq I_{n-r} \oplus \diag \left( \psi_r,\ldots, \psi_1 \right) \in \FR^{n\times n}$.
\end{description}
\end{thm}
\smallskip \noindent
\textbf{Proof}: The proof is based on the claim that there exist two right coprime matrices $\widetilde{A} \in \FR^{n\times n}$, with $\det \widetilde{A} \ne 0$, and $\widetilde{C} \in \FR^{p\times n}$ with Smith forms, respectively, $S_A$ and $S_C$ as in the statement and such that $R = \widetilde{C} \widetilde{A}^{-1}\in \F^{p\times n}$. Then, {Proposition} \ref{thm.1stpraux1}-ii) implies that $A$ and $\widetilde{A}$ have the same Smith form and that $C$ and $\widetilde{C}$ have also the same Smith form, which completes the argument.

It remains to prove the claim. Let $U\in \FR^{p\times p}$ and $V\in \FR^{n\times n}$ be unimodular matrices such that $R = U S_R V$ and let us define
\begin{align*}
  \widetilde{C} & := U \, \left( \diag \left(\varepsilon_1,\ldots, \varepsilon_r \right) \oplus 0_{(p-r) \times (n-r)} \right) \in \FR^{p\times n},  \\
  \widetilde{A} & := V^{-1} \, \left( \diag \left( \psi_1,\ldots, \psi_r \right) \oplus  I_{n-r} \right) \in \FR^{n\times n}.
\end{align*}
Obviously $\widetilde{C}$ has Smith form $S_C$, $\widetilde{A}$ has Smith form $S_A$, $\det \widetilde{A} \ne 0$ and $R = \widetilde{C} \widetilde{A}^{-1}$. It remains to prove that $\widetilde{C}$ and $\widetilde{A}$ are right coprime. To this purpose, note that $\gcd(\varepsilon_i,\psi_i) \doteq 1$ for all $i=1,\dots,r$. Moreover, since $\FR$ is a B\'{e}zout domain, there exist $x_i, y_i\in \FR$ such that $x_i \varepsilon_i + y_i \psi_i = 1$, for $i=1,\ldots r$. At this point, define the matrices
\begin{align*}
   X & := \left( \diag \left(x_1,\ldots, x_r \right) \oplus 0_{(n-r) \times (p-r)} \right) \, U^{-1} \in \FR^{n\times p},  \\
   Y & :=  \left( \diag \left( y_1,\ldots, y_r \right) \oplus  I_{n-r} \right) \, V \in \FR^{n\times n},
\end{align*}
which obviously satisfy $X \widetilde{C} + Y \widetilde{A} = I_n$.
\qed \smallskip

A result analogous to Theorem \ref{thm.1stpraux2} can be given for $R = A^{-1} B$ if $A$ and $B$ are left coprime.

We now prove Theorem \ref{thm.1stpraux3}. Observe that the existence of the matrix $U$ in its statement is guaranteed by the left coprime version of Proposition \ref{teo_copr_an}-iii). Theorem \ref{thm.1stpraux3} is related to \cite[Lemma 6.3-8]{Kail80} (stated and proved for univariate polynomial and rational matrices with complex coefficients). However, we emphasize that Theorem \ref{thm.1stpraux3} goes beyond \cite[Lemma 6.3-8]{Kail80} not only because of its weaker assumptions on the base ring, but also for extending the statement to include some results that are essential for proving Theorem \ref{thm.1stproof}.

\begin{thm}  \label{thm.1stpraux3}
Let $A\in \FR^{n\times n}$, with $\det A\neq 0$, and  $B\in \FR^{n\times m}$ be left coprime and let
\be\label{eq.UQfirst}
U=\begin{bmatrix} Y_{11} & Y_{12}\\Y_{21} & Y_{22}\end{bmatrix} \in\FR^{(n+m)\times (n+m)}
\ee
be a unimodular matrix such that $\begin{bmatrix}A & B\end{bmatrix} U=\begin{bmatrix} I_n & 0\end{bmatrix}$. Let
$$
S_{A^{-1} B}  \doteq \diag \left(\frac{\gamma_1}{\theta_1},\ldots, \frac{\gamma_s}{\theta_s}\right) \oplus 0_{(n-s) \times (m-s)} \in \F^{n\times m}
$$
be the Smith-McMillan form of $A^{-1} B \in \F^{n \times m}$. Then:
\begin{description}
  \item[\rm i)] $\det Y_{22} \ne 0$.
  \item[\rm ii)] The Smith form of $A$ is $S_A \doteq I_{n-s} \oplus \diag \left( \theta_s,\ldots, \theta_1 \right)\in \FR^{n\times n}$.
  \item[\rm iii)] The Smith form of $Y_{22}$ is $S_{Y_{22}} \doteq I_{m-s} \oplus \diag \left( \theta_s,\ldots, \theta_1 \right)\in \FR^{m\times m}$.
  \item[\rm iv)] The Smith form of $B$ is $S_B \doteq \diag \left(\gamma_1,\ldots, \gamma_s\right) \oplus 0_{(n-s) \times (m-s)} \in \FR^{n\times m}$.
  \item[\rm v)] The Smith form of $Y_{12}$ is $S_{Y_{12}} \doteq \diag \left(\gamma_1,\ldots, \gamma_s\right) \oplus 0_{(n-s) \times (m-s)} \in \FR^{n\times m}$.
\end{description}
Thus, the nontrivial invariant factors of $A$ and $Y_{22}$ coincide,  while $B$ and $Y_{12}$ are equivalent.
\end{thm}

\smallskip \noindent
\textbf{Proof}:
\noindent [\,i)\,] The equality  $\begin{bmatrix}A & B\end{bmatrix} U=\begin{bmatrix} I_n & 0\end{bmatrix}$ implies
\begin{equation}\label{eq.1-aux3}
A Y_{12} + B Y_{22} = 0.
\end{equation}
Thus, any $x\in \F^{m}$ such that $Y_{22} x = 0$ satisfies $A Y_{12} x= 0$, and hence $Y_{12} x= 0$  because $\det A \ne 0$. Thus
$U \left[ \begin{smallmatrix}
0 \\
x
\end{smallmatrix}\right] =0$, which implies that $x =0$ because $U$ is unimodular. Therefore, $Y_{22}$ is nonsingular, i.e., $\det Y_{22} \ne 0$.

\bigskip \noindent
[\,ii) and iv)\,] These items are corollaries of the left coprime version of Theorem \ref{thm.1stpraux2} applied to $A^{-1} B$.

\bigskip \noindent
[\,iii) and v)\,] From \eqref{eq.1-aux3} we see that $Y_{12} Y_{22}^{-1} = - A^{-1} B$, and hence $Y_{12}Y_{22}^{-1}$ has the same Smith-McMillan form as $A^{-1}B$.
Moreover, Proposition \ref{teo_copr_an}-iv) and the fact that $U$ is unimodular imply that $Y_{12}$ and $Y_{22}$ are right coprime. Thus, items iii) and v) follow from applying Theorem \ref{thm.1stpraux2} to $Y_{12} Y_{22}^{-1}$. \qed \medskip

We are finally in the position of proving Rosenbrock's theorem over EDDs.
As mentioned in the introduction, Rosenbrock's theorem relates the invariant factors of given matrices $A\in\FR^{n\times n}$ and
$P=\begin{bsmallmatrix}A & B\\C& D\end{bsmallmatrix}\in\FR^{(n+p)\times(n+m)}$ with the invariant fractions of the Schur
complement of $A$ in $P$, that is, $G=D-CA^{-1}B$, when $A$ and $B$ are left coprime and, moreover, $A$ and $C$ are right
coprime. It is important to bear in mind that if $r=\rank G$ then it may be that either $n\geq r$ or $r\geq n$
(or both, i.e, $r=n$).  For example,
if $A$ is unimodular and  $D=CA^{-1}B$ (note that, in this case, $D\in\FR^{p\times m}$) then $r=0$ and so  $n>r$.
On the other hand, it could be $\min(m,p)>n$ and $r>n$. A simple example
is: $D=I_2$, $A=1$, $B=\begin{bsmallmatrix}1 & 0\end{bsmallmatrix}$ and $C=\begin{bsmallmatrix}0\\1
\end{bsmallmatrix}$. In this case $n=1$ and $r=2$. Finally, we can take $A=1$, $B=C=0$ and $D=1$ to construct an example where $n=r=1$.

\begin{rem}
{\rm To follow the statement of Theorem \ref{thm.1stproof}, as well as the analysis in Section \ref{sec.beyond} below, it is useful to keep in mind that the matrix $P$ in \eqref{eq.firstdefP} always satisfies
$\rank P = n + \rank G$, independently of the coprimeness properties of $A,B$ and $A,C$. This is a consequence of the
nonsingularity of $A$ and the block LU factorization (over $\F$)
\begin{equation} \label{eq.rankPG}
\begin{bmatrix}
  A & B \\
  C & D
\end{bmatrix} =
\begin{bmatrix}
  I_n & 0 \\
  CA^{-1} & I_p
\end{bmatrix}\begin{bmatrix}
  A & 0 \\
  0 & G
\end{bmatrix}\begin{bmatrix}
  I_n & A^{-1} B \\
  0 & I_m
\end{bmatrix}.
\end{equation}  }
\end{rem}

\begin{thm} {\rm (Rosenbrock's theorem over EDDs)} \label{thm.1stproof}
Let $A\in\FR^{n\times n}$, $B\in\FR^{n\times m}$, $C\in\FR^{p\times n}$ and $D\in\FR^{p\times m}$ with $\det A\neq 0$. Let
\begin{equation} \label{eq.firstdefP}
P=\begin{bmatrix} A & B\\ C & D\end{bmatrix}\in\FR^{(n+p)\times (n+m)}
\end{equation}
and
\begin{equation} \label{eq.firstdefG}
G=D-CA^{-1}B\in\F^{p\times m}, \qquad r=\rank G.
\end{equation}
Assume that $A$ and $B$ are left coprime and that $A$ and $C$ are right coprime.
\begin{enumerate}
\item If the Smith-McMillan form of $G$ is
$$
S_G \doteq \diag \left(\frac{\varepsilon_1}{\psi_1},\ldots, \frac{\varepsilon_r}{\psi_r}\right) \oplus 0_{(p-r) \times (m-r)} \in \F^{p\times m},
$$
and $g$ is the largest index in $\{1,\ldots ,r \}$ such that $\psi_g  \notin {\FR^\times}$, then
%and $g=\min(n,r)$, then
\begin{description}
  \item[\rm i)] the Smith form of $P$ is $$S_P \doteq I_n \oplus \diag \left(\varepsilon_1,\ldots, \varepsilon_r \right) \oplus 0_{(p-r) \times (m-r)} \in \FR^{(n+p)\times (n+m)}, \quad \mbox{and}$$
  \item[\rm ii)] the Smith form of $A$ is $S_A \doteq I_{n-g} \oplus \diag \left( \psi_g,\ldots, \psi_1 \right) \in \FR^{n\times n}$.
% \item[\rm iii)]  if $g=n$ then $\psi_{n+1}\eqdot\cdots\eqdot\psi_r\eqdot 1$.
\end{description}
\item Conversely, if the Smith forms of $A$ and $P$ are $S_A\doteq\diag \left( \psi_n,\ldots, \psi_1 \right) \in \FR^{n\times n}$ and
$S_P\doteq\diag \left(\varepsilon_1,\ldots, \varepsilon_{n+r} \right)\oplus 0_{(p-r) \times (m-r)} \in \FR^{(n+p)\times (n+m)}$,  respectively, then
\begin{description}
\item[\rm iii)]  $\varepsilon_1\eqdot \cdots\eqdot\varepsilon_n\eqdot1$,
\item[\rm iv)] $\psi_n\eqdot\cdots\eqdot\psi_{r+1}\eqdot 1$  if $n\geq r$, and
\item[\rm v)] the Smith-McMillan form of $G$ is $$S_G\eqdot  \diag \left(\frac{\varepsilon_{n+1}}{\psi_1},\ldots,
\frac{\varepsilon_{n+r}}{\psi_r}\right) \oplus 0_{(p-r) \times (m-r)} \in \F^{p\times m}$$ where,  if $r\geq n$,
$\psi_{n+1}\eqdot\cdots\eqdot\psi_r\eqdot 1$ .
 \end{description}
 \end{enumerate}
\end{thm}

\smallskip
\textbf{Proof}: Let
$
U\in\FR^{(n+m)\times (n+m)}
$,
partitioned as in \eqref{eq.UQfirst},
be a unimodular matrix such that $\begin{bmatrix}A & B\end{bmatrix} U=\begin{bmatrix} I_n & 0\end{bmatrix}$. Then
\[
PU=\begin{bmatrix} I_n & 0\\CY_{11}+DY_{21} & CY_{12}+DY_{22}\end{bmatrix}
\]
and
\begin{equation}\label{eq.1-firstprrosen}
\begin{bmatrix} I_n & 0\\-CY_{11}-DY_{21} & I_p\end{bmatrix} PU=
\begin{bmatrix} I_n & 0\\0 & CY_{12}+DY_{22}\end{bmatrix}.
\end{equation}
Hence, $P$ and $I_n \oplus (CY_{12}+DY_{22})$ are equivalent. The rest of the proof consists of relating $G \in \F^{p\times m}$ and $CY_{12}+DY_{22}\in \FR^{p\times m}$. Note first that $\begin{bmatrix}A & B\end{bmatrix} U=\begin{bmatrix} I_n & 0\end{bmatrix}$ implies \eqref{eq.1-aux3}, which in turn implies
$G Y_{22}  = D Y_{22} + C Y_{12}$. Moreover, $Y_{22}$ is nonsingular by Theorem \ref{thm.1stpraux3}. Therefore,
\begin{equation}\label{eq.2-firstprrosen}
G = (CY_{12}+DY_{22}) \, Y_{22}^{-1}.
\end{equation}
The next step is to prove that $(CY_{12}+DY_{22})$ and $Y_{22}$ are right coprime. For this purpose note that $Y_{12}$ and $Y_{22}$ are right coprime by Proposition \ref{teo_copr_an}-iv). Therefore, there exist $E_1 \in \FR^{m\times n}$ and $E_2 \in \FR^{m\times m}$ such that $I_m = E_1 Y_{12} + E_2 Y_{22}$. Moreover, since $A$ and $C$ are right coprime, there exist $F_1 \in \FR^{n\times n}$ and $F_2 \in \FR^{n\times p}$ such that $I_n = F_1 A + F_2 C$. Therefore,
\begin{align*}
  I_m & =  E_1 Y_{12} + E_2 Y_{22} \\
      & =  E_1 ( F_1 A + F_2 C) Y_{12} + E_2 Y_{22}\\
      & =  E_1 F_2 C Y_{12} + (E_2 -E_1 F_1B )  Y_{22}\\
      & =  E_1 F_2 \, (C Y_{12} + D Y_{22}) + (E_2 -E_1 F_1B - E_1 F_2 D )  Y_{22} ,
\end{align*}
where the third equality follows from \eqref{eq.1-aux3}. Thus, $(CY_{12}+DY_{22})$ and $Y_{22}$ are right coprime
by Proposition \ref{teo_copr_an}-v). Then, item i) follows from applying Theorem \ref{thm.1stpraux2} to \eqref{eq.2-firstprrosen} and combining the result with \eqref{eq.1-firstprrosen}. Theorem \ref{thm.1stpraux2} applied to \eqref{eq.2-firstprrosen} also implies that  the Smith form of $Y_{22}$ is $I_{m-r}\oplus\diag\left(\psi_r,\ldots,\psi_1\right)=
I_{m-g}\oplus\diag\left(\psi_g,\ldots,\psi_1\right)$ (note that $m\geq r$ and
by definition $\psi_{g+1}\doteq \cdots\doteq\psi_r\doteq 1$).
On the other hand, according to items ii) and  iii) of Theorem \ref{thm.1stpraux3}, the Smith forms of $A$  and $Y_{22}$
differ only in the number of trivial invariant factors. Hence,
$S_A \doteq I_{n-g} \oplus \diag \left( \psi_g,\ldots, \psi_1 \right)$.

Conversely, items iii), iv) and v) follow from items i) and ii) and from the uniqueness of the Smith forms of $A$ and $P$
up to a product by  diagonal unimodular matrices.\qed\smallskip

\section{Beyond Rosenbrock's theorem} \label{sec.beyond} In this section we investigate the relations between the Smith-McMillan form of $G$ in \eqref{eq.firstdefG} and the Smith forms of $A$ and $P$ in \eqref{eq.firstdefP}, when $A$ and $B$ are not left coprime or $A$ and $C$ are not right coprime, i.e., when  the coprimeness assumptions of Rosenbrock's theorem do not hold. We will see that in this case item ii) in Theorem \ref{thm.1stproof} never holds and that the same happens for item i) under some additional conditions which are important in applications.

{ To study these questions, we need two additional auxiliary results. Proposition \ref{pro.Lem1Co} is Lemma 1 of  \cite{Coppel74}. Although \cite{Coppel74} focuses on PIDs, the proof of \cite[Lemma 1]{Coppel74} only requires basic properties of the invariant factors that can be easily seen to hold over any EDD, namely, the existence of a Smith form for each $A\in \FR^{p\times m}$, the definition of the determinantal divisors and the fact that the $k$th determinantal divisor is the product of the $k$th ``smallest''
invariant factors up to a product by a unit of $\FR$, i.e., equation \eqref{eq.ifdetdiv}. Hence, it holds for matrices over EDDs. Proposition \ref{pro.Lem1Co} also follows from the stronger result \cite[Theorem 4.3]{Shchedrykbook}, that provides necessary and sufficient conditions for a matrix with entries in an EDD to be a left divisor of another matrix. More information on problems related to the invariant factors of products of matrices can be found in \cite{CaQue15} and \cite[Chapter 4]{Shchedrykbook}.}

\begin{propos}\label{pro.Lem1Co}
Let $A_1\in\FR^{m\times n}$, $A_2\in\FR^{n\times p}$ and let $A=A_1A_2$. Let $\alpha_1^{(1)}\mid\cdots\mid
\alpha_{r_1}^{(1)}$, $\alpha_1^{(2)}\mid\cdots\mid \alpha_{r_2}^{(2)}$ and $\alpha_1\mid\cdots\mid
\alpha_{r}$ be the invariant factors of $A_1$, $A_2$ and $A$, respectively. Then $\alpha_k^{(j)}\mid \alpha_k$
for $j=1,2$ and $k=1,\ldots, r$.
\end{propos}

Since Rosenbrock's theorem deals with Schur complements \cite[{Section 0.8.5}]{H-J13}, the (known) Theorem \ref{thm.schurcomplement} below is a natural tool in this context. It is stated, assuming that the system matrix $P$ has entries in a field, in \cite[p. 25]{H-J13}. However, it is generally valid for matrices with entries in any integral domain since an integral domain can be embedded in its
field of fractions. Theorem \ref{thm.schurcomplement} uses the following notation borrowed from \cite[Section 1.7]{Fried16}. For a positive integer $n$ denote $[n] = \{1,2,\ldots , n\}$ and for $k \in [n]$ denote by $[n]_k$ the set of all subsets of $[n]$ having cardinality $k$. Each element $\mathI \in [n]_k$ is represented as $\mathI = (i_1 , \ldots , i_k)$ with
$1\leq i_1 < \cdots < i_k \leq n$. For any matrix $H$ of size $s \times t$ and for any $\mathI \in [s]_k$,
$\mathJ \in [t]_\ell$, we denote by $H[\mathI, \mathJ]$ the $k\times \ell$ matrix with entries that lie in the rows
of $H$ indexed by $\mathI$ and in the columns indexed by $\mathJ$.

\begin{thm} \label{thm.schurcomplement} Let $\FR$ be any integral domain. Let $A\in\FR^{n\times n}$, $B\in\FR^{n\times m}$, $C\in\FR^{p\times n}$ and $D\in\FR^{p\times m}$ with $\det A\neq 0$. Let
\[
P=\begin{bmatrix} A & B\\ C & D\end{bmatrix}\in\FR^{(n+p)\times (n+m)}
\quad
\mbox{and}
\quad
G=D-CA^{-1}B\in\F^{p\times m}.
\]
Let $\mathI = (i_1 , \ldots , i_k)  \in [p]_k$, $\mathJ = (j_1 , \ldots , j_k) \in [m]_k$ and $\widehat{\mathI} = (i_1 +n , \ldots , i_k +n)$, $\widehat{\mathJ} = (j_1 + n, \ldots , j_k + n)$. Then
\begin{equation} \label{eq.schur1}
\det P[\, [n] \cup \widehat{\mathI} , [n] \cup \widehat{\mathJ} \, ] = (\det A) \, (\det G [\mathI, \mathJ]).
\end{equation}
Moreover, for any $e \in \FR$,
\begin{equation} \label{eq.schur2} (
\det (e P)[\, [n] \cup \widehat{\mathI} , [n] \cup \widehat{\mathJ} \,]) = (\det (e A)) \,
(\det (e G) [\mathI, \mathJ]).
\end{equation}
\end{thm}

\medskip

The results in this section are based on Theorem \ref{thm.beyond1} below, which is a generalization of \cite[Theorem 7]{Coppel74}. The proof in \cite{Coppel74} is valid over UFDs but not over EDDs that are not PIDs. Theorem \ref{thm.beyond1} goes beyond \cite[Theorem 7]{Coppel74} not only for its different assumptions on the base ring,
but also for its detailed analysis of the properties of the matrices $E$ and $F$: this is crucial for other results in this section.

\begin{thm} \label{thm.beyond1}
Let $A\in\FR^{n\times n}$, $B\in\FR^{n\times m}$, $C\in\FR^{p\times n}$ and $D\in\FR^{p\times m}$ with $\det A\neq 0$. If $A$ and $B$ are not left coprime or $A$ and $C$ are not right coprime, %{\color{purple} (or both)}
then there exist matrices $A_0\in\FR^{n\times n}$,  with $\det A_0 \neq 0$, $B_0\in\FR^{n\times m}$, $C_0 \in\FR^{p\times n}$, $E\in\FR^{n\times n}$ and $F\in \FR^{n\times n}$ such that
\begin{equation}\label{eq.from0tone0}
\begin{bmatrix} A & B\\ C & D\end{bmatrix}=
\begin{bmatrix} E & 0\\ 0 & I_p\end{bmatrix}
\begin{bmatrix} A_0 & B_0\\ C_0 & D\end{bmatrix}
\begin{bmatrix} F & 0\\ 0 & I_m \end{bmatrix}
\end{equation}
and
\begin{description}
  \item[\rm i)] $A_0$ and $B_0$ are left coprime and $A_0$ and $C_0$ are right coprime;
  \item[\rm ii)] $\det E \ne 0$, $\det F \ne 0$, and at least one of these determinants is not a unit of $\FR$;
    \item[\rm iii)] $D-CA^{-1}B = D-C_0 A_0^{-1}B_0$.
    \end{description}
Moreover, if a factorization of the form \eqref{eq.from0tone0} satisfies properties i), ii), iii) above, then:
  \begin{description}
  \item[\rm iv)] if $A$ and $B$ are not left coprime then $E$ can be chosen with $\det E \notin {\FR^\times}$;
  \item[\rm v)] if $A$ and $C$ are not right coprime then $F$ can be chosen with $\det F \notin {\FR^\times}$;
  \item[\rm vi)] if $A$ and $B$ are left coprime then $E$ must be unimodular and $E$ can be chosen to be $I_n$;
  \item[\rm vii)] if $A$ and $C$ are right coprime then $F$ must be unimodular and $F$ can be chosen to be $I_n$.

\end{description}
\end{thm}

\medskip
\textbf{Proof}. We first prove items i), ii), iii) and iv), splitting the proof in two steps. The first one constructs the matrix $E$ and the second one the matrix $F$.

Step 1. If $A$ and $B$ are left coprime, we set $E=I_n$. If $A$ and $B$ are not left coprime, then every Smith form
$\begin{bmatrix} S & 0\end{bmatrix} \in \FR^{n\times (n+m)}$ of $\begin{bmatrix} A & B\end{bmatrix}$ satisfies $\det S \ne 0$, because $\det A \ne 0$, and $\det S \notin {\FR^\times}$, because of Proposition \ref{teo_copr_an}-ii). Hence we can write
\begin{equation} \label{eq.1from0to}
\begin{bmatrix} A & B\end{bmatrix} = U \begin{bmatrix} S & 0\end{bmatrix}
\begin{bmatrix}
  W_{11} & B_0 \\
  W_{21} & W_{22}
\end{bmatrix} = U S
\begin{bmatrix}
  W_{11} & B_0
\end{bmatrix},
\end{equation}
where $U$ and $\begin{bsmallmatrix}
  W_{11} & B_0 \\
  W_{21} & W_{22}
\end{bsmallmatrix} \in \FR^{(n+m) \times (n+m)}$ are unimodular. In this case we take $E = U S$, which satisfies $\det E \notin {\FR^\times}$, because $\det S \notin {\FR^\times}$, and $\det E \ne 0$. Note that the definition of $E$ and \eqref{eq.1from0to} allow us to write
\begin{equation} \label{eq.2from0to}
\begin{bmatrix} A & B\\ C & D\end{bmatrix}=
\begin{bmatrix} E & 0\\ 0 & I_p\end{bmatrix}
\begin{bmatrix} W_{11} & B_0 \\ C & D\end{bmatrix},
\end{equation}
with $W_{11}$ and $B_0$ left coprime by Proposition \ref{teo_copr_an}-iv) and $\det W_{11} \ne 0$ because $A = E W_{11}$.

\medskip

Step 2. If $W_{11}$ and $C$ in \eqref{eq.2from0to} are right coprime, we set $F = I_n$ and \eqref{eq.2from0to} provides a factorization as in \eqref{eq.from0tone0} with properties i), ii), iii) and iv).\footnote{We observe \emph{en passant} that it is possible that $W_{11}$ and $C$ are right coprime and, at the same time, $A$ and $C$ are not right coprime: see Remark \ref{rem.noconverse} and Example \ref{ex:noconverse}.}

\smallskip

If $W_{11}$ and $C$ are not right coprime, then every Smith form $\begin{bsmallmatrix} T \\ 0\end{bsmallmatrix} \in \FR^{(n + p)\times n}$ of $\begin{bsmallmatrix} W_{11} \\ C \end{bsmallmatrix}$ satisfies $\det T \ne 0$, because $\det W_{11} \ne 0$, and $\det T \notin {\FR^\times}$, because of Proposition \ref{teo_copr_an}-ii). Hence we can write
\begin{equation} \label{eq.3from0to}
\begin{bmatrix} W_{11} \\ C \end{bmatrix} = \begin{bmatrix}
  A_0 & X_{12} \\
  C_0 & X_{22}
\end{bmatrix} \begin{bmatrix} T \\ 0\end{bmatrix}
V =
\begin{bmatrix}
  A_0 \\ C_0
\end{bmatrix} TV,
\end{equation}
where $V$ and $\begin{bsmallmatrix}
  A_0 & X_{12} \\
  C_0 & X_{22}
\end{bsmallmatrix} \in \FR^{(n+p) \times (n+p)}$ are unimodular. In this case we take $F = TV$, which satisfies $\det F \notin {\FR^\times}$, because $\det T \notin {\FR^\times}$, and $\det F \ne 0$. Finally, the definition of $F$, \eqref{eq.2from0to} and \eqref{eq.3from0to} allow us to write
\begin{equation} \label{eq.4from0to}
\begin{bmatrix} A & B\\ C & D\end{bmatrix}=
\begin{bmatrix} E & 0\\ 0 & I_p\end{bmatrix}
\begin{bmatrix} A_0 & B_0 \\ C_0 & D\end{bmatrix}
\begin{bmatrix} F & 0\\ 0 & I_m\end{bmatrix},
\end{equation}
with $A_0$ and $C_0$ right coprime by Proposition \ref{teo_copr_an}-iv). Moreover, $A_0$ and $B_0$ are left coprime, because each common left divisor of $A_0$ and $B_0$ is a common left divisor of $W_{11} = A_0 F$ and $B_0$, and these two matrices are left coprime. Then, \eqref{eq.4from0to} provides a factorization  as in \eqref{eq.from0tone0} satisfying properties i), ii), iii) and iv).

In order to prove item v), we can repeat the same argument but switching the order of steps 1 and 2.

It remains to prove items vi) and vii). The argument above already shows that if $A$ and $B$ are left coprime then it is possible to choose $E=I_n$ and that if $A$ and $C$ are right coprime then it is possible to choose $F=I_n$. (For the latter statement, observe that if $A$ and $C$ are right coprime, then $W_{11}$ and $C$ are right coprime because each common right divisor of $W_{11}$ and $C$ is a common right divisor of $A=EW_{11}$ and $C$.) To show that $E$ must be unimodular when $A$ and $B$ are left coprime, suppose that \eqref{eq.from0tone0} holds with $\det E \not\in {\FR^\times}$; then $A=E(A_0F)$ and $B=EB_0$, implying that $E$ is a non-unimodular common left divisor of $A$ and $B$, thus contradicting their left coprimeness. This completes the proof of item vi); the proof of item vii) can be finalized in an analogous manner from the equalities $A=(EA_0)F$ and $C=C_0F$.
\qed\smallskip

\begin{rem} \label{rem.noconverse}
{\rm The properties of $E$ in Theorem \ref{thm.beyond1}-iv) and of $F$ in Theorem \ref{thm.beyond1}-v) are not guaranteed to hold simultaneously. In other words, even if $A$ and $B$ are not left coprime \emph{and} $A$ and $C$ are not right coprime, it is not guaranteed that $E$ and $F$ can be chosen with $\det E \notin {\FR^\times}$ \emph{and} $\det F \notin {\FR^\times}$. Instead, in this scenario, Theorem \ref{thm.beyond1} only implies the weaker consequence that \emph{at least one} between $E$ and $F$ can be picked to not be unimodular. We illustrate this with Example \ref{ex:noconverse} below.}
\end{rem}

\begin{ex}\label{ex:noconverse}
{\rm In this example, we show that, if $\FR$ contains at least one prime\footnote{There exist elementary divisor domains that do not contain any prime, e.g., the algebraic integers \cite[Remark 3.10]{N23}.}, then one may be unable to pick $E,F$ \emph{both} non-unimodular even when neither $(A,B)$  nor $(A,C)$ are coprime.

Consider the following matrices:
\[
\begin{bmatrix}
  A & B \\
  C & D
\end{bmatrix} =
\left[ \begin{array}{cc|cc}
         p & 0 & p & 0 \\
         0 & 1 & 0 & 1 \\ \hline
         p & 1 & 0 & 0
       \end{array}
\right].
\]
Observe that $\begin{bsmallmatrix} p & 0 \\ 0 & 1\end{bsmallmatrix}$ is a common left divisor of $A$ and $B$ and a common right divisor of $A$ and $C$ with $\det \begin{bsmallmatrix} p & 0 \\ 0 & 1\end{bsmallmatrix} = p$. Hence, if $p  \notin {\FR^\times}$, then $A$ and $B$ are not left coprime and $A$ and $C$ are not right coprime.

Suppose now that $p$ is prime, and take any factorization of the form \eqref{eq.from0tone0} with $\det E \not \in {\FR^\times} \not \ni \det F$ and the pairs $(A_0,B_0)$ and $(A_0,C_0)$ both coprime. This is impossible because then, $(\det E \det F)$ divides $\det A = p$, which implies that at most one between $\det E$ and $\det F$ can be a non-unit. 

Instead, if $p$ is neither a unit nor a prime, then it \emph{is} possible to pick $E,F$ both non-unimodular in  \eqref{eq.from0tone0}. Namely, write $p=ab$ with $a \not \in {\FR^\times} \not \ni b$ and pick
\[  E = \begin{bmatrix}
a&0\\
0&1
\end{bmatrix}, F=\begin{bmatrix}
b &0\\
0&1
\end{bmatrix} \Rightarrow  \begin{bmatrix}
A & B\\
C & D
\end{bmatrix} = (E \oplus 1)  \left[ \begin{array}{cc|cc}
         1 & 0 & b & 0 \\
         0 & 1 & 0 & 1 \\ \hline
         a & 1 & 0 & 0
       \end{array} \right] (F \oplus I_2)   .\]}
\end{ex}

Let us now extend to an arbitrary EDD the following nomenclature that was originally introduced for $\FR=\C[z]$ in the linear systems theory literature \cite{Coppel74,Kail80,Rosen70}. A representation of a matrix $G \in \F^{p\times m}$ in terms of matrices over $\FR$ as $G = D - C A^{-1} B$ is called a {\em realization} of $G$. If $A$ and $B$ are left coprime and $A$ and $C$ are right coprime, then the representation is called an {\em irreducible realization} (or a realization of least order or a minimal realization). As a corollary of Theorem \ref{thm.beyond1} we obtain that every matrix with entries in the field of fractions of an EDD has an irreducible realization.

\begin{corol} \label{cor.existencereal} Every matrix $G \in \F^{p\times m}$ can be written as $G = D - C A^{-1} B$ with $A\in\FR^{n\times n}$, $B\in\FR^{n\times m}$, $C\in\FR^{p\times n}$ and $D\in\FR^{p\times m}$, with $\det A\neq 0$, with $A$ and $B$ left coprime and with $A$ and $C$ right coprime.
\end{corol}
\medskip
\textbf{Proof}. If $G \in \FR^{p\times m}$, then $D = 0$, $C = -G$, $A= I_m$, and $B=I_m$ trivially satisfy the conditions of the statement. Otherwise, let $\varphi$ be the least common denominator of the entries of $G$ and observe that $\varphi \notin {\FR^\times}$. Then, $G = (\varphi G) (\varphi I_m)^{-1}$. This expression is of the form $G = D - C A^{-1} B$ with $D = 0$, $C = -\varphi G$, $A= \varphi I_m$, and $B=I_m$, but $A$ and $C$ are not right coprime. However, Theorem \ref{thm.beyond1} can be applied to these matrices $A,B,C,D$ to obtain other matrices $A_0,B_0,C_0,D$ with the properties in the statement and $G = D - C_0 A_0^{-1} B_0$.
\qed\smallskip

The idea of the proof of Corollary \ref{cor.existencereal} is taken from that in \cite[pp. 96-97]{Coppel74}, which is only valid when $\FR$ is a PID. Another proof of Corollary \ref{cor.existencereal} valid for an arbitrary EDD $\FR$ is given by the construction of
$\widetilde{A}$ and $\widetilde{C}$ in the proof of Theorem \ref{thm.1stpraux2}.

Next, we present the main result of this section. To this goal, let us formally define $\psi_i:=1$ for $r+1 \leq i \leq n$ (note that this definition is only relevant when $r < n$, as the condition on the index $i$ is empty if $r \geq n$. If $r<n$, this definition extends the divisibility chain of the $\psi_i$, i.e., now it holds that $\psi_{i+1} \mid \psi_i$ for all $i=1,\dots,n-1$.)
\begin{thm} \label{thm.mainbeyond}
Let $A\in\FR^{n\times n}$, $B\in\FR^{n\times m}$, $C\in\FR^{p\times n}$ and $D\in\FR^{p\times m}$ with $\det A\neq 0$,
\begin{equation} \label{eq.Pmainbeyond}
P=\begin{bmatrix} A & B\\ C & D\end{bmatrix}\in\FR^{(n+p)\times (n+m)}, \quad \mbox{and} \quad
G=D-CA^{-1}B\in\F^{p\times m}.
\end{equation}
Let
$$
S_G \doteq \diag \left(\frac{\varepsilon_1}{\psi_1},\ldots, \frac{\varepsilon_r}{\psi_r}\right) \oplus 0_{(p-r) \times (m-r)} \in \F^{p\times m},
$$
be the Smith-McMillan form of $G$ and 
$g$ be the largest index in $\{1,\ldots ,r \}$ such that $\psi_g \notin {\FR^\times}$. Assume that $A$ and $B$ are not left coprime or that $A$ and $C$ are not right coprime. If the Smith form of $A$ is
\[
S_A \doteq \diag \left( \widetilde{\psi}_n \, , \ldots , \widetilde{\psi}_1  \right)
\]
and the Smith form of $P$ is
\[
S_P \doteq \diag \left( \widetilde{\varepsilon}_1 , \ldots , \widetilde{\varepsilon}_{n+r}  \right) \oplus 0_{(p-r) \times (m-r)}
\]
then
\begin{description}
  \item[\rm i)] $n \geq  g$ and $\psi_i \mid \widetilde{\psi}_i$, for $i= 1, \ldots, n$;
  \item[\rm ii)] $\displaystyle \frac{\widetilde{\psi}_n  \cdots  \widetilde{\psi}_2 \, \widetilde{\psi}_1 }{\psi_g  \cdots \psi_2 \, \psi_1 } \notin {\FR^\times}$, that is,  $\frac{\widetilde{\psi}_i}{\psi_i} \notin {\FR^\times}$ for at least one value of $i=1,\dots,n$;
  \item[\rm iii)] $\varepsilon_i \mid \widetilde{\varepsilon}_{n+i}$ for $i=1, \ldots , r$;
  \item[\rm iv)] $\displaystyle
      \frac{\widetilde{\varepsilon}_1  \widetilde{\varepsilon}_2 \cdots \widetilde{\varepsilon}_{n+r}}{\varepsilon_1 \varepsilon_2 \cdots \varepsilon_r }
      \mid
      \frac{\widetilde{\psi}_n  \cdots  \widetilde{\psi}_2 \, \widetilde{\psi}_1 }{\psi_g  \cdots \psi_2 \, \psi_1 };
      $
  \item[\rm v)] if $r=p=m$, i.e., if $G$ and $P$ are square and nonsingular, then
      $$
      \frac{\widetilde{\varepsilon}_1  \widetilde{\varepsilon}_2 \cdots \widetilde{\varepsilon}_{n+r}}{\varepsilon_1 \varepsilon_2 \cdots \varepsilon_{r}}
      \doteq
      \frac{\widetilde{\psi}_n  \cdots  \widetilde{\psi}_2 \, \widetilde{\psi}_1 }{\psi_g  \cdots \psi_2 \, \psi_1 } \notin {\FR^\times},
      $$
       that is, $\widetilde{\varepsilon}_i \notin {\FR^\times}$ for at least one $i = 1, \ldots , n$, or $\displaystyle \frac{\widetilde{\varepsilon}_{n+i}}{\varepsilon_i}\notin {\FR^\times}$ for at least one $i = 1, \ldots , r$.

\end{description}
\end{thm}
\textbf{Proof}. Theorem \ref{thm.beyond1} guarantees that $P$ can be factorized as in \eqref{eq.from0tone0}, where the matrices $A_0, B_0, C_0, E$ and $F$ have the properties established in Theorem \ref{thm.beyond1}. Moreover, Theorem \ref{thm.1stproof} applied to $P_0= \begin{bsmallmatrix} A_0 & B_0 \\ C_0 & D \end{bsmallmatrix}$ and Theorem \ref{thm.beyond1}-iii) imply that the Smith forms of $P_0$ and $A_0$ are, respectively,
\begin{align*}
  S_{P_0}  &  \doteq I_n \oplus \diag \left(\varepsilon_1,\ldots, \varepsilon_r \right) \oplus 0_{(p-r) \times (m-r)}, \\
  S_{A_0} & \doteq I_{n-g} \oplus \diag \left( \psi_g,\ldots, \psi_1 \right).
\end{align*}
With these ingredients at hand we proceed to prove all the items in the statement.

\medskip \noindent i) The expression of $S_{A_0}$ implies that $n\geq g$.
 The factorization \eqref{eq.from0tone0} implies $A = E A_0 F$. Combining this equality with Proposition \ref{pro.Lem1Co} and the expressions of $S_{A_0}$ and $S_{A}$ yields the result.

\medskip \noindent ii) $A = E A_0 F$ and the expressions of $S_{A_0}$ and $S_{A}$ imply
\begin{equation}\label{eq.aux1mainbeyond}
\widetilde{\psi}_n  \cdots  \widetilde{\psi}_2 \, \widetilde{\psi}_1 \doteq \det A = (\det E) (\det A_0 ) (\det F) \doteq (\psi_g  \cdots \psi_2 \, \psi_1)  (\det E) (\det F),
\end{equation}
and the result follows from the fact that at least one of $\det E$ and $\det F$ is not a unit of $\FR$.

\medskip \noindent iii) This item follows from combining \eqref{eq.from0tone0} with Proposition \ref{pro.Lem1Co} and the expressions of $S_{P_0}$ and $S_{P}$.

\medskip \noindent iv) To prove this item, we apply \eqref{eq.schur2} with $e = \psi_1$ and $k=r$. For this purpose note first that $\psi_1 G \in \FR^{p\times m}$ and has Smith form $$
S_{\psi_1 G } \doteq \diag \left(\psi_1 \, \frac{\varepsilon_1}{\psi_1},\ldots, \psi_1 \, \frac{\varepsilon_r}{\psi_r}\right) \oplus 0_{(p-r) \times (m-r)}.
$$ Moreover, the Smith form of $\psi_1 P$ is
$$S_{\psi_1 P} \doteq \diag \left(\psi_1  \widetilde{\varepsilon}_1 , \ldots , \psi_1 \widetilde{\varepsilon}_{n+r}  \right) \oplus 0_{(p-r) \times (m-r)}.$$
Then, by \eqref{eq.ifdetdiv}, the greatest common divisor of all the minors of size $(n+r) \times (n+r)$ of $\psi_1 P$ is $\psi_1^{n+r} \widetilde{\varepsilon}_1  \cdots \widetilde{\varepsilon}_{n+r}$. %, up to a product by a unit.
 Therefore, from \eqref{eq.schur2} we get
\begin{equation}\label{eq.beyond2iv}
\psi_1^{n+r} \widetilde{\varepsilon}_1  \cdots \widetilde{\varepsilon}_{n+r} \mid  (\det (\psi_1 A)) \,
(\det (\psi_1 G) [\mathI, \mathJ])
\end{equation}
for all $\mathI \in [p]_r$ and $\mathJ \in [m]_r$. Note that
\begin{align*}
\gcd & \left\{ (\det (\psi_1 A))
(\det (\psi_1 G) [\mathI, \mathJ]) \; : \; \mathI \in [p]_r ,  \mathJ \in [m]_r\right\}\\ & \phantom{aaaaaaaa}  \doteq (\det (\psi_1 A)) \gcd \left\{
\det (\psi_1 G) [\mathI, \mathJ]  \; : \; \mathI \in [p]_r\, , \, \mathJ \in [m]_r\right\} \\
& \phantom{aaaaaaaa} \doteq (\det (\psi_1 A)) \, \psi_1^r \, \frac{\varepsilon_1}{\psi_1} \cdots \frac{\varepsilon_r}{\psi_r} \\ & \phantom{aaaaaaaa}  \doteq
\psi_1^{n+r} \, \frac{\widetilde{\psi}_n  \cdots  \widetilde{\psi}_2 \, \widetilde{\psi}_1 }{\psi_g  \cdots \psi_2 \, \psi_1 } \, \varepsilon_1 \cdots \varepsilon_r,
\end{align*}
where the first ``equality up to a product by a unit'' follows from Lemma \ref{lem.div1}, the second one from \eqref{eq.ifdetdiv} and the third one from the expression of $S_A$. Combining the equality above with \eqref{eq.beyond2iv} yields the result.

\medskip \noindent v) If $P$ is nonsingular, then $P_0= \begin{bsmallmatrix} A_0 & B_0 \\ C_0 & D \end{bsmallmatrix}$ is also nonsingular by \eqref{eq.from0tone0} and $\det P =  (\det E) (\det P_0 ) (\det F)$. Combining this with \eqref{eq.aux1mainbeyond} yields $$\frac{\det P}{\det P_0} = \frac{\det A}{\det A_0} =  (\det E) (\det F)\notin {\FR^\times}.$$ The result follows from the expressions of the Smith forms $S_P, S_{P_0}, S_A$ and $S_{A_0}$
\qed\smallskip

\begin{rem} \label{rem.applbeyond} {\rm
Note that given $P_0= \begin{bsmallmatrix} A_0 & B_0 \\ C_0 & D \end{bsmallmatrix}$, with $A_0$ and $B_0$ left coprime and $A_0$ and $C_0$ right coprime, equation \eqref{eq.from0tone0} with arbitrary nonsingular matrices $E,F \in \FR^{n\times n}$ may yield a matrix $P = \begin{bsmallmatrix} A & B \\ C & D \end{bsmallmatrix}$ such that the Smith forms of $A$ and $P$ are completely different from the denominators and numerators of the Smith-McMillan form of the Schur complement $G$ of $A$ in $P$, according to Theorem \ref{thm.mainbeyond}-ii) and v) and their proofs (as mentioned in the introduction to Proposition \ref{pro.Lem1Co}, see \cite{CaQue15} for a thorough analysis of the  relationship between the invariant factors of a product of matrices over EDDs and of its factors).
In any case the divisibility relations in Theorem \ref{thm.mainbeyond}-i) and iii) are always guaranteed.  Suppose that $\FR=\C[z]$ is the ring of univariate polynomials with complex coefficients and $\F=\C(z)$ is the field of corresponding rational functions. Then, potentially, this remark has practical (negative) implications for algorithms \cite{nleigs} that compute the zeros of a rational matrix $G$ via a linear polynomial system matrix $P$ whose Schur complement (or transfer function matrix) is $G$ without guaranteeing that the coprimeness conditions of the pairs $(A,B)$ and $(A,C)$ hold.}
\end{rem}

It is natural to wonder whether items iii) and iv) in Theorem \ref{thm.mainbeyond} can be strengthened to guarantee that at least one of the involved divisions is proper, i.e., at least one of the involved quotients is not a unit. In fact, this is the case under the additional nonsingularity assumption of item v). The next example illustrates that Theorem \ref{thm.mainbeyond} cannot be strengthened in that direction and that there are matrices $P$ with $A,B$ not left coprime and $A,C$ not right coprime whose nontrivial invariant factors are exactly equal to the numerators of the invariant fractions of $G$.

\begin{ex} {\rm Let $\varepsilon_1, \psi_1, \alpha_2, \alpha_3 \in \FR$ be nonzero and be not units. Assume also that  $\gcd (\varepsilon_1, \psi_1) \eqdot 1$. Consider
$$
P = \begin{bmatrix}
      A & B \\
      C & D
    \end{bmatrix} =
    \left[ \begin{array}{ccc|ccc}
             -\psi_1 & 0 & 0 & 1 & 0 & 0 \\
             0 & \alpha_2 & 0 & 0 & 1 & 0 \\
             0 & 0 & \alpha_3 & 0 & 0 & 0 \\ \hline
             \varepsilon_1 & 0 & 0 & 0 & 0 & 0 \\
             0 & 0 & 0 & 0 & 0 & 0 \\
             0 & 0 & 1 & 0 & 0 & 0
           \end{array}
    \right].
$$
Note that
$$
G = D - C A^{-1} B = \begin{bmatrix}
                       \frac{\varepsilon_1}{\psi_1} & 0 & 0 \\
                       0 & 0 & 0 \\
                       0 & 0 & 0
                     \end{bmatrix}
$$ is directly the Smith-McMillan form.
Applying elementary row and column operations, it is easy to find that the Smith form of $P$ is $S_P \doteq I_3 \oplus \varepsilon_1 \oplus 0_{2\times 2}$. Moreover,
$$
\begin{bmatrix}
1 & 0 & 0 \\
0 & 1 & 0 \\
0 & 0 & \alpha_3
\end{bmatrix}
\quad
\mbox{and}
\quad
\begin{bmatrix}
1 & 0 & 0 \\
0 & \alpha_2 & 0 \\
0 & 0 & 1
\end{bmatrix}
$$
are, respectively, a common non-unimodular left divisor of $A$ and $B$ and
a common non-unimodular right divisor of $A$ and $C$.
}
\end{ex}

Under additional full rank assumptions on $G$, we can complete Theorem \ref{thm.mainbeyond}-iv) with Theorem \ref{thm.bismainbeyond}. Note that the results in Theorem \ref{thm.bismainbeyond} are not as strong as Theorem \ref{thm.mainbeyond}-v), because we cannot guarantee the equality of the quotient of the $\varepsilon_i$'s and of the quotient of the $\psi_i$'s.

\begin{thm} \label{thm.bismainbeyond}
Let $A\in\FR^{n\times n}$, $B\in\FR^{n\times m}$, $C\in\FR^{p\times n}$ and $D\in\FR^{p\times m}$ with $\det A\neq 0$,
\[
P=\begin{bmatrix} A & B\\ C & D\end{bmatrix}\in\FR^{(n+p)\times (n+m)}, \quad \mbox{and} \quad
G=D-CA^{-1}B\in\F^{p\times m}.
\]
Let
$$
S_G \doteq \diag \left(\frac{\varepsilon_1}{\psi_1},\ldots, \frac{\varepsilon_r}{\psi_r}\right) \oplus 0_{(p-r) \times (m-r)} \in \F^{p\times m}
$$
be the Smith-McMillan form of $G$ and
\[
S_P \doteq \diag \left( \widetilde{\varepsilon}_1 , \ldots , \widetilde{\varepsilon}_{n+r}  \right) \oplus 0_{(p-r) \times (m-r)}
\]
be the Smith form of $P$.
\begin{description}
  \item[\rm i)] If $A$ and $B$ are not left coprime and $r = \rank G = p$, that is, $G$ and $P$ have full row rank, then
  $\displaystyle \frac{\widetilde{\varepsilon}_1  \widetilde{\varepsilon}_2 \cdots \widetilde{\varepsilon}_{n+p}}{\varepsilon_1 \varepsilon_2 \cdots \varepsilon_p} \notin {\FR^\times}$.
  \item[\rm ii)] If $A$ and $C$ are not right coprime and $r = \rank G = m$, that is, $G$ and $P$ have full column rank, then $\displaystyle \frac{\widetilde{\varepsilon}_1  \widetilde{\varepsilon}_2 \cdots \widetilde{\varepsilon}_{n+m}}{\varepsilon_1 \varepsilon_2 \cdots \varepsilon_m} \notin {\FR^\times}$.
\end{description}
\end{thm}
\textbf{Proof}.
We only prove item i), because the proof of item ii) then follows by transposing every relevant matrix.  Theorem \ref{thm.beyond1} guarantees that $P$ can be factorized as in \eqref{eq.from0tone0}, with $\det E \notin {\FR^\times}$, $A_0, B_0$ left coprime and  $A_0,C_0$ right coprime. Moreover, Theorem \ref{thm.1stproof} applied to $P_0= \begin{bsmallmatrix} A_0 & B_0 \\ C_0 & D \end{bsmallmatrix}$ and Theorem \ref{thm.beyond1}-iii) imply that the Smith form of $P_0$ is
\begin{align*}
  S_{P_0}  &  \doteq I_n \oplus \diag \left(\varepsilon_1,\ldots, \varepsilon_p \right) \oplus 0_{0 \times (m-p)}.
\end{align*}
From the factorization \eqref{eq.from0tone0} and the Cauchy-Binet formula we see that each minor of order $(n+p)$ of $P$ is the product of $\det E$ times the corresponding minor of $P_0 \, (F \oplus I_m)$.
Thus, using Lemma \ref{lem.div1}, we get the following relation between determinantal divisors
\[
(\det E) \, \delta_{n+p} (P_0 \, (F \oplus I_m)) \eqdot\delta_{n+p} (P).
\]
Moreover, it follows from the Cauchy-Binet identity that $\delta_{n+p} (P_0) \mid \delta_{n+p} (P_0 (F \oplus I_m))$. Therefore, $(\det E) \,\delta_{n+p} (P_0) \mid \delta_{n+p} (P)$. From \eqref{eq.ifdetdiv} and the expressions of $S_{P_0}$ and $S_P$, we get
$$
(\det E) \, (\varepsilon_1 \varepsilon_2 \cdots \varepsilon_p) \mid (\widetilde{\varepsilon}_1  \widetilde{\varepsilon}_2 \cdots \widetilde{\varepsilon}_{n+p}),
$$
and the result follows.
\qed\smallskip

We conclude this section by completing the first part of Rosenbrock's theorem to a
necessary and sufficient condition.  It is also possible to complete the second part of Rosenbrock's theorem to a necessary and sufficient condition. We omit the statement of that result for brevity.

\begin{thm}  \label{thm.strongrosen}
Let $A\in\FR^{n\times n}$, $B\in\FR^{n\times m}$, $C\in\FR^{p\times n}$ and $D\in\FR^{p\times m}$ with $\det A\neq 0$. Let
\[
P=\begin{bmatrix} A & B\\ C & D\end{bmatrix}\in\FR^{(n+p)\times (n+m)},
\qquad
G=D-CA^{-1}B\in\F^{p\times m},
\]
and
$$
S_G \doteq \diag \left(\frac{\varepsilon_1}{\psi_1},\ldots, \frac{\varepsilon_r}{\psi_r}\right) \oplus 0_{(p-r) \times (m-r)} \in \F^{p\times m}
$$
be the Smith-McMillan form of $G$, where 
$g$ is the largest index in $\{1,\ldots ,r \}$ such that $\psi_g \notin {\FR^\times}$. Consider the conditions
\begin{description}
  \item[\rm a)] $A$ and $B$ are left coprime,
  \item[\rm b)] $A$ and $C$ are right coprime,
\end{description}
and
\begin{description}
  \item[\rm i)] the Smith form of $P$ is $S_P \doteq I_n \oplus \diag \left(\varepsilon_1,\ldots, \varepsilon_r \right) \oplus 0_{(p-r) \times (m-r)} \in \FR^{(n+p)\times (n+m)},$
  \item[\rm ii)] the Smith form of $A$ is $S_A \doteq I_{n-g} \oplus \diag \left( \psi_g,\ldots, \psi_1 \right) \in \FR^{n\times n}$.
\end{description}
Then, {\rm a)} and {\rm b)} both hold if and only if {\rm i)} and {\rm ii)} both hold.
\end{thm}
\textbf{Proof}. Theorem \ref{thm.1stproof} proved that ``{\rm a)} and {\rm b)}''  implies ``{\rm i)} and {\rm ii)}''.  Theorem \ref{thm.mainbeyond} proved that not ``{\rm a)} and {\rm b)}''  implies not ``{\rm ii)}", and hence, not ``{\rm i)} and {\rm ii)}''.
\qed\smallskip

\vspace*{-10pt}%
\subsection{A further improvement of Theorem \ref{thm.mainbeyond} over PIDs} If $\FR$ is a PID, then it is simultaneously an EDD and a UFD. Therefore any matrix $A\in \FR^{p\times m}$ has a Smith form as in \eqref{eq.smith} and each of its invariant factors $\alpha_i$ has a unique (up to order) factorization as a finite product of prime elements of $\FR$. This allows us to write
\begin{align}
  \alpha_1 & = \beta_1^{e_{11}} \beta_2^{e_{12}} \cdots \beta_\ell^{e_{1\ell}}, \nonumber \\
  \alpha_2 & = \beta_1^{e_{21}} \beta_2^{e_{22}} \cdots \beta_\ell^{e_{2\ell}}, \nonumber \\
  \vdots &  \phantom{aaaaaaaaaaa} \vdots \label{eq.elemdivisors} \\
  \alpha_r & = \beta_1^{e_{r1}} \beta_2^{e_{r2}} \cdots \beta_\ell^{e_{r\ell}}, \nonumber
\end{align}
where $\beta_1, \ldots, \beta_\ell$ are prime elements of $\FR$ and $e_{ij}$ are nonnegative integers that satisfy
$$
0 \leq e_{1j} \leq e_{2j} \leq \cdots \leq e_{rj}, \quad j=1,\ldots , \ell.
$$

The factors $\beta_j^{e_{ij}}$ with $e_{ij}>0$ are called the elementary divisors of $A$ \cite[{Definition 5, Chapter VI}]{Gantmacher59}. Given any prime element $\pi \in \FR$, we will consider the finite sequence of the associated \emph{partial multiplicities} of $A$ at $\pi$, i.e., the sequence of the positive integers $t_i$ such that $\alpha_i = \pi^{t_i} \, \gamma_{i}$ with  $\gamma_{i} \in \FR$, and $\gcd (\pi, \gamma_i) \doteq 1$ for $i=1,\ldots, r$.
We emphasize that this sequence is empty when {$\pi \nmid  \alpha_i$, $i = 1,\ldots , r$} %$\pi \not\doteq \beta_j$, $j = 1,\ldots , \ell$;
otherwise, the sequence is non-empty and non-decreasing (since $\alpha_i$ divides $\alpha_{i+1}$ for $i=1, \ldots , r-1$). Observe that if all the elementary divisors of $A\in \FR^{p\times m}$ are known and the {\em rank} of $A$ is also known, then the invariant factors of $A$ can be constructed via \eqref{eq.elemdivisors}. Therefore the knowledge of the elementary divisors or, equivalently, of their prime factors and the associated partial multiplicities of a matrix over a PID is essentially equivalent to the knowledge of its invariant factors.

If $\F$ is the field of fractions of the PID $\FR$ and $G\in \F^{p\times m}$ has the Smith-McMillan form \eqref{eq.smithmcmillan}, then the elementary divisors of $\diag \left(\varepsilon_1,\ldots, \varepsilon_r \right)$ will be called the {\em numerator elementary divisors} of $G$ and the elementary divisors of
$\diag \left(\psi_1,\ldots, \psi_r\right)$ will be called the {\em denominator elementary divisors} of $G$. These elementary divisors are very important and well-known when $\FR$ is the ring of univariate polynomials with coefficients in a field and $\F$ is the field of corresponding rational functions \cite{AmMaZa14,Kail80,Rosen70}. In this setting, they are sometimes called zero and pole elementary divisors, respectively \cite{local-dmqvd}.

With these definitions at hand, we prove the following corollary of Theorem \ref{thm.mainbeyond} where
the following property will be used: if the Smith-McMillan form of $G=D-CA^{-1}B$ is
$\diag \left(\frac{\varepsilon_1}{\psi_1},\ldots, \frac{\varepsilon_r}{\psi_r}\right) \oplus 0_{(p-r) \times (m-r)}$ then
$\psi_1\cdots\psi_g\mid \det A$. This is a consequence of Theorem \ref{thm.mainbeyond} -i) and the fact that $\det A$ is the
product of the invariant factors of $A$ up to a product by a unit.

\begin{corol} \label{cor.mainbeyond} Let $\FR$ be a PID and $\F$ its field of fractions. Let $A\in\FR^{n\times n}$, $B\in\FR^{n\times m}$, $C\in\FR^{p\times n}$ and $D\in\FR^{p\times m}$ with $\det A\neq 0$,
\[
P=\begin{bmatrix} A & B\\ C & D\end{bmatrix}\in\FR^{(n+p)\times (n+m)}, \quad \mbox{and} \quad
G=D-CA^{-1}B\in\F^{p\times m}.
\]
Let
$$
S_G \doteq \diag \left(\frac{\varepsilon_1}{\psi_1},\ldots, \frac{\varepsilon_r}{\psi_r}\right) \oplus 0_{(p-r) \times (m-r)} \in \F^{p\times m}
$$
be the Smith-McMillan form of $G$ and $g$ be the largest index in $\{1,\ldots ,r \}$ such that $\psi_g \notin {\FR^\times}$.
If $\pi \in\FR$ is prime and
$\displaystyle \gcd \left( \pi \, ,\,  \frac{\det A}{\psi_g  \cdots \psi_2 \, \psi_1 } \right) \doteq 1$, then the sequence of the partial multiplicities of $P$ at $\pi$ is equal to the sequence of the partial multiplicities of the numerators of the invariant fractions of $G$ at $\pi$.
\end{corol}
\textbf{Proof}. It suffices to prove the statement when either $A$ and $C$ are not right coprime or $A$ and $B$ are not left coprime (or both), as otherwise the statement follows immediately by Rosenbrock's theorem. We will use the notation as in the statement of Theorem \ref{thm.mainbeyond}. Let $\widetilde{\varepsilon}_i = \pi^{t_i} \, \gamma_{i}$, with $t_i \geq 0$ integer, $\gamma_i \in \FR$, $\gcd (\pi, \gamma_i) \doteq 1$, for $i = 1, \ldots, n+r$, and let $\varepsilon_i = \pi^{s_i} \, \phi_{i}$, with $s_i \geq 0$ integer, $\phi_i \in \FR$, $\gcd (\pi, \phi_i) \doteq 1$, for $i = 1, \ldots, r$. Theorem \ref{thm.mainbeyond}-iii) implies $s_i \leq t_{n+i}$ for $i=1,\ldots , r$. Moreover, combining items iii) and iv) of Theorem \ref{thm.mainbeyond}, we get
\[
\varepsilon_1 \varepsilon_2 \cdots \varepsilon_r \mid
\widetilde{\varepsilon}_1  \widetilde{\varepsilon}_2 \cdots \widetilde{\varepsilon}_{n+r}
      \mid
      \frac{\det A}{\psi_g  \cdots \psi_2 \, \psi_1 } \, \varepsilon_1 \varepsilon_2 \cdots \varepsilon_r.
\]
Since  $\displaystyle \gcd \left( \pi \, ,\,  \frac{\det A}{\psi_g  \cdots \psi_2 \, \psi_1 } \right) \doteq 1$, we see that
$s_1 + \cdots + s_r = t_1 + \cdots + t_{n+r}$. This equality and $s_i \leq t_{n+i}$ for $i=1,\ldots , r$ imply $t_1 = \cdots = t_n =0$ and $s_i = t_{n+i}$ for $i=1,\ldots , r$, which gives the result.
\qed\smallskip

Observe that the result in Corollary \ref{cor.mainbeyond} can be applied to any prime element $\pi$ that satisfies $\gcd \left( \pi \, ,\, \det A \right) \doteq 1$. This has practical applications in situations where $\det A$ has just a few prime factors \cite{lietaertetal2022,subai11} that can be easily identified, because in this case $P$ may contain most of the ``zero'' structure of $G$ even when $A$ and $B$ are not left coprime or $A$ and $C$ are not right coprime. The ``local'' (at $\pi$) result in Corollary \ref{cor.mainbeyond} is in the spirit of the local results presented in \cite{local-dmqvd} for polynomial system matrices of rational matrices. A wealth of information on ``local'' Smith-McMillan forms of rational matrices can be found in \cite{AmMaZa14}, though not applied specifically to polynomial system matrices.

\section{Three corollaries of Rosenbrock's theorem} \label{sec.consequences}
In this section we present three consequences of Theorem \ref{thm.1stproof} related to the denominators of the Smith-McMillan form of matrices in $\F^{p\times m}$. To express these results concisely, we adopt the convention that if $G\in\F^{p\times m}$, $r=\rank G$ and $\psi_r\mid\psi_{r-1}\mid\cdots\mid \psi_1$
are the denominators of its invariant fractions then $\psi_k=1$ for $k > r$.

The first result in this section establishes that the 
denominators of the invariant fractions of a matrix with entries in $\F$ do not change by the addition of any matrix  with entries in $\FR$.

\begin{corol}\label{cor.Th1Co}
Let $G_1\in\F^{p\times m}$, $P_2\in\FR^{p\times m}$ and $G=G_1+P_2$. Let
$\psi_{r_1}^{(1)}\mid\psi_{r_1-1}^{(1)}\mid\cdots\mid \psi_1^{(1)}$ and
$\psi_r\mid\psi_{r-1}\mid\cdots\mid \psi_1$ be the denominators of the invariant fractions of $G_1$ and $G$,
respectively. 
Then, for $k\geq 1$, $\psi_k \doteq \psi_k^{(1)}$.
\end{corol}
\textbf{Proof}. Corollary \ref{cor.existencereal} implies that $G_1$ can be written as $G_1 = D - C A^{-1} B$ with $A\in\FR^{n\times n}$, $B\in\FR^{n\times m}$, $C\in\FR^{p\times n}$ and $D\in\FR^{p\times m}$, with $\det A\neq 0$, with $A$ and $B$ left coprime and with $A$ and $C$ right coprime. Then $G$ can be written as $G = (D+P_2) - C A^{-1} B$ with $D+P_2\in\FR^{p\times m}$. The result follows from applying Theorem \ref{thm.1stproof}-ii) first to $P_1 = \begin{bsmallmatrix} A & B \\ C & D \end{bsmallmatrix}$ and $G_1$ and then to
$P = \begin{bsmallmatrix} A & B \\ C & D +P_2 \end{bsmallmatrix}$ and $G$.
\qed\smallskip

\hide{
\begin{rem}\label{rem.CopThm11}{\rm
Recall that $G\in\F^{p\times m}$ always admits realizations.  That is to say,  $G$ can be written as $G=D-CA^{-1}B$ for some
matrices $A$, $B$, $C$ and $D$ with entries in $\FR$ (see Corollary \ref{cor.existencereal}
and its introduction). Following \cite{Coppel74}, if $A\in\FR^{n\times n}$, we say that  $n$ is
the \textit{dimension} of the realization. Theorem 11 of \cite{Coppel74} says that when $\FR$ is a PID then
the quantity $g$ defined in the statement of Theorem \ref{thm.1stproof} is the \textbf{minimal dimension }of any realization
of $G$. One can show  that using Corollary \ref{cor.Th1Co}  the proof of that theorem still holds when
the underlying ring $\FR$ is an EDD.
}
\end{rem}
}

The second result in this section is similar to the first one but deals with matrices and submatrices. This result is used without proof in \cite[pp. 114-115]{Vergh79} when $\FR$ is the ring of univariate polynomials with real coefficients as an intermediate step when generalizing Rosenbrock's theorem to capture the structure at infinity via strongly irreducible polynomial system matrices.
\begin{corol}\label{lem.Verghese}
Let $T\in \FR^{n\times s}$, $U\in \FR^{n\times m}$, $V\in \FR^{p\times s}$, $W\in\F^{p\times m}$
and
\[
R=\left[ \begin{array}{rr} T  & U\\ V & W\end{array} \right] \in \F^{(n+p) \times (s+m)}.
\]
Let $r=\rank W$ and $l=\rank R$. Let $\displaystyle \frac{\varepsilon_1}
{\psi_1}$, \ldots, $\displaystyle \frac{\varepsilon_r}{\psi_r}$
be the invariant fractions of $W$ and let $\displaystyle \frac{\epsilon_1}{\phi_1}$, \ldots,
$\displaystyle \frac{\epsilon_l}{\phi_l}$ be the invariant fractions of $R$. Then $l\geq r$ and
$\psi_j \doteq \phi_j$ for $ j\geq 1$.
\end{corol}
\textbf{Proof}.  It is plain that $l\geq r$. Now, observe that
\[
R= \left[ \begin{array}{rc} T & U\\V & 0\end{array} \right]+\begin{bmatrix} 0 & 0\\0 & W\end{bmatrix},
\]
where $\begin{bsmallmatrix}T & U\\V & 0\end{bsmallmatrix}$ and
$\begin{bsmallmatrix} 0 & 0\\0 & W\end{bsmallmatrix}$ have their entries in $\FR$ and $\F$,
respectively. It follows from Corollary \ref{cor.Th1Co} that $\psi_j \doteq \phi_j$ for $j\geq 1$.
\qed\smallskip

Finally, recall that any $G\in\F^{p\times m}$ admits realizations.  That is to say,  $G$ can be written as $G=D-CA^{-1}B$ for some
matrices $A$, $B$, $C$ and $D$ with entries in $\FR$ (see Corollary \ref{cor.existencereal}
and its introduction). Following \cite{Coppel74}, if $A\in\FR^{n\times n}$, we say that  $n$ is
the \textit{dimension} of the realization. Theorem 11 of \cite{Coppel74} says that when $\FR$ is a PID then
the quantity $g$ defined in the statement of Theorem \ref{thm.1stproof} is the \textit{minimal dimension} of any realization
of $G$. We prove in the next corollary that this result remains true when
the underlying ring $\FR$ is an EDD.

\begin{corol} \label{cor.g} Let $G \in \F^{p \times m}$, $\frac{\varepsilon_1}{\psi_1} , \ldots , \frac{\varepsilon_r}{\psi_r}$ be its invariant fractions and $g$ be the largest index in $\{1, \ldots , r \}$ such that
$\psi_g \notin {\FR^\times}$. Then $g$  is the minimal dimension of any realization of $G$.
\end{corol}
\textbf{Proof}. Let $G = D - C A^{-1} B$ be any realization of $G$ with $A \in \FR^{n\times n}, B \in \FR^{n\times m}, C \in \FR^{p\times n}, D \in \FR^{p\times m}$ and $\det A \ne 0$. If $A,B$ are left coprime and $A,C$ are right coprime, then Theorem \ref{thm.1stproof}-ii) implies $n \geq g$. Otherwise, Theorem \ref{thm.mainbeyond}-i) implies $n \geq g$. It only remains to find a realization of $G$ of dimension $g$. For that, we follow the construction in the proof of \cite[Theorem 11]{Coppel74}. Let
$$
U^{-1} G V^{-1} = \diag \left(\frac{\varepsilon_1}{\psi_1},\ldots, \frac{\varepsilon_r}{\psi_r}\right) \oplus 0_{(p-r) \times (m-r)}
$$
be the Smith-McMillan form of $G$, where $U$ and $V$ are unimodular. Then
\begin{align*}
  A & = \diag \left(\psi_1, \ldots , \psi_g \right) \in \FR^{g\times g}, \\
  B  & = \begin{bmatrix}
          I_g & 0_{g \times (m-g)}
        \end{bmatrix} V \in \FR^{g\times m},\\
  C & = U \begin{bmatrix}
            \diag \left(\varepsilon_{1},\ldots, \varepsilon_g \right) \\
            0_{(p-g) \times g}
          \end{bmatrix} \in \FR^{p \times g}, \\
  D & = U \, ( 0_{g\times g} \oplus \diag \left(\varepsilon_{g+1},\ldots, \varepsilon_r \right) \oplus 0_{(p-r) \times (m-r)}) \, V \in \FR^{p\times m},
\end{align*}
where $D =0$ if $g=r$, gives a realization $G = D - C A^{-1} B$ of $G$ of dimension $g$. Moreover, $A,B$ are left coprime and $A,C$ are right coprime, since $\gcd ( \varepsilon_i, \psi_i) \doteq 1$ for $i = 1, \ldots , r$ (see the proof of Theorem \ref{thm.1stpraux2}).
\qed\smallskip

\section{Rosenbrock's theorem for matrices over $\F$}\label{sec.RTFF}
In this section $G\in\F^{p\times m}$ will be a matrix written in the form $G=D-CA^{-1}B$ where $A$, $B$, $C$ and $D$
are themselves matrices with entries in $\F$. Our goal is to investigate whether, and to what extent Rosenbrock's theorem still holds
in this setting. The motivation of this investigation is the following one: Assume that $\K$ is an arbitrary field,
$\K[z]$ the ring of polynomials in the unknown $z$ and $\K(z)$ the field of rational functions. Let
$P(z)=\begin{bsmallmatrix}A(z) & B(z)\\C(z) & D(z)\end{bsmallmatrix}\in\K[z]^{(n+p)\times(n+m)}$ be
an irreducible polynomial system matrix with  $A(z)\in\K[z]^{n\times n}$ and $\det A(z)\neq 0$. Then
$G(z)=D(z)-C(z)A(z)^{-1}B(z)$ is the transfer function matrix of $P(z)$ and we aim to compute the poles and zeros at infinity
of $G(z)$. The natural ring where the Smith-McMillan form of $G(z)$ exhibits its poles and zeros at infinity is the ring
of proper rational functions (see \cite{AmMaZa14,AmMaZa16,DiCo82,NV23laa}  and \cite[Chapter 3]{Vardu91} for more references),  i.e.
\[
\K_{pr}(z)=\left\{\frac{p(z)}{q(z)}:p(z),q(z)\in\K[z], \deg(q(z))\geq \deg(p(z))\right\}.
\]
The zeros and poles at infinity of matrix polynomials in the Smith-McMillan sense are closely related to their eigenvalues at infinity in the Gohberg-Lancaster-Rodman sense, see for example
\cite{AmMaZa14,AnDoHoMa23,NV23laa}. On the other hand, when $\K \subseteq \C$, there are efficient and stable algorithms to compute the eigenvalues at
infinity  of matrix polynomials, for instance, using an  appropriate linearization and then solving the  corresponding  generalized
eigenvalue problem via the QZ algorithm \cite{moler-stewart}. However, it is plain
that nonconstant  matrix  polynomials do not have their entries in $\K_{pr}(z)$ but in its field of fractions $\K(z)$.  This is a
reason why we think that  a ``Rosenbrock-like
theorem" when the matrices $A$, $B$, $C$ and $D$ have their entries in the field of fractions of $\FR$ is of interest.
Another, more theoretical in nature, reason to study  this problem is that, to our knowledge, such extensions of Rosenbrock's theorem remain fully unexplored. 

An important assumption in Rosenbrock's theorem is the coprimeness of $A$, $B$ and $A$, $C$.
How can this assumption be adapted to matrices over a field of fractions? A possible answer, suggested in \cite{AmMaZa16}, is given in Definition \ref{rl_cop}. For brevity, in the rest of the paper, we call the least common denominator of the entries of a matrix $H \in \F^{p \times m}$ simply the least common denominator of $H$.

\begin{defi}\label{rl_cop}
Let $\FR$ be an EDD and $\F$ its field of fractions. Let $T_1\in\F^{p\times m}$ and $~T_2\in\F^{q\times m}$
and let $\tau_1, \tau_2$ be  least common denominators in $\FR$ of $T_1$ and $T_2$,
respectively. $T_1$ and $T_2$ are said to be right coprime in $\FR$ if  $\tau_1T_1 \in \FR^{p\times m}$ and $\tau_2T_2 \in \FR^{q\times m}$ are right
coprime in $\FR$. Similarly, $T_1\in\F^{m\times p}$ and $~T_2\in\F^{m\times q}$ are left coprime in $\FR$ if
 $\tau_1T_1 \in \FR^{m\times p}$ and $\tau_2T_2 \in \FR^{m\times q}$ are left coprime in $\FR$.
\end{defi}

Observe that Definition \ref{rl_cop} emphasizes that $\tau_1, \tau_2$ are least common denominators {\em in the ring $\FR$}. The reason is that different integral domains may be different subrings of the same field of fractions. Thus a matrix with entries in this field of fractions
may have different least common denominators in those different rings.  For example,  $\K(z)$ happens to be the field of fractions of both
$\K_{pr}(z)$ and $\K[z]$. Consider the  rational matrix
\[
T(z)=\begin{bmatrix} \frac{z^2+1}{(z-1)^2(z-2)^2} & \frac{z+3}{(z-1)^3}\\
\frac{z^3}{(z-2)^4(z-1)} & \frac{z+5}{(z-1)^3(z-2)}\end{bmatrix}\in\C(z)^{2\times 2}.
\]
A least common denominator of $T(z)$ in  $\C_{pr}(z)$ is $1$ because this matrix is in  ${\C_{pr}(z)^{2 \times 2}}$. However, a least common
denominator  in  $\C[z]$ is $(z-1)^3(z-2)^4$, which is not an associate of $1$ in $\C[z]$.

It can be shown \cite{AmMaZa16} that if $T(z)\in\K[z]^{m\times n}$ and $d_T$ is its degree then a least common
denominator of $T(z)$  in $\K_{pr} (z)$ is  $\left(\frac{1}{z}\right)^{d_T}$. Indeed, $\frac{1}{z^{d_T}}T(z)\in\K_{pr}(z)^{m\times n}$.

\begin{thm} {\rm (First part of Rosenbrock's theorem over fields of fractions of EDDs)}
\label{thm.1stproof-1}  Let $\FR$ be an EDD and $\F$ its field of fractions.
Let $A\in\F^{n\times n}$, $B\in\F^{n\times m}$ and $C\in\F^{p\times n}$
be matrices such that $\det A\neq 0$, $A$ and $B$ are left coprime in $\FR$ and $A$ and $C$ are right coprime in $\FR$.
Let $\alpha$, $\beta$ and $\gamma$ be least common denominators in $\FR$ of $A$,
$B$ and $C$, respectively.  Let  $D\in\F^{p\times m}$ be a matrix such that
$\frac{\beta\gamma}{\alpha}D\in\FR^{p\times m}$. Let
$G=D-CA^{-1}B\in\F^{p\times m}$
and
$$
S_G \doteq \diag \left(\frac{\varepsilon_1}{\psi_1},\ldots, \frac{\varepsilon_r}{\psi_r}\right) \oplus 0_{(p-r) \times (m-r)}
\in \F^{p\times m},
$$
be its Smith-McMillan form.  Then
\begin{description}
\item[\rm i)] the Smith-McMillan form of $A$ is $S_A \doteq \frac{1}{\alpha}I_{(n-g)} \oplus
\diag \left( \frac{\wt{\psi}_g}{\delta_g\chi_g},\ldots, \frac{\wt\psi_1}{\delta_1\chi_1} \right) \in \F^{n\times n}$, and
\item[\rm ii)] the Smith-McMillan form of
\begin{equation} \label{eq.firstdefP-1}
P=\frac{1}{\beta\gamma}\begin{bmatrix} \alpha A & \beta B\\ \gamma C & \frac{\beta\gamma}{\alpha} D\end{bmatrix}
\in \F^{(n+p)\times (n+m)}
\end{equation}
is $S_P\doteq\frac{1}{\beta\gamma}I_n \oplus \diag \left(\frac{\wt\varepsilon_1}{\nu_1\chi_1},\ldots,
\frac{\wt\varepsilon_r}{\nu_r\chi_r} \right) \oplus 0_{(p-r) \times (m-r)} \in \F^{(n+p)\times (n+m)}$,
\end{description}
where
\be\label{eq.deltai}
\delta_i=\gcd(\varepsilon_i,\alpha), \quad \wh\alpha_i=\frac{\alpha}{\delta_i}, \quad \wt\varepsilon_i=\frac{\varepsilon_i}{\delta_i},
\quad i=1,\ldots, r,
\ee
\be\label{eq.nui}
\nu_i=\gcd(\psi_i,\beta\gamma),\quad \wh\beta_i=\frac{\beta\gamma}{\nu_i},\quad \wt{\psi}_i=\frac{\psi_i}{\nu_i},
\quad i=1,\ldots, r,
\ee
\be\label{eq.chii}
\chi_i=\gcd(\wh\alpha_i,\wh\beta_i),\quad  \wt\alpha_i=\frac{\wh\alpha_i}{\chi_i},\quad \wt\beta_i=\frac{\wh\beta_i}{\chi_i},
\quad i=1,\ldots, r,
\ee
and $g$ is the largest index in $\{1,\ldots, r\}$ such that $\wt\alpha_g\wt\psi_g\notin{\FR^\times}$.
\end{thm}

\textbf{Proof}. Let $Q=\beta\gamma P=\begin{bsmallmatrix} \alpha A & \beta B\\ \gamma C &
\frac{\beta\gamma}{\alpha}D\end{bsmallmatrix} \in\FR^{(n+p)\times (n+m)}$ and
$\wt{G}=\frac{\beta\gamma}{\alpha}D- \gamma C(\alpha A)^{-1}\beta B=\frac{\beta\gamma}{\alpha}(D-CA^{-1}B)=
\frac{\beta\gamma}{\alpha}G$. Then, there are unimodular matrices $U\in\FR^{p\times p}$ and $V\in\FR^{m\times m}$
such that
\[
U\wt{G}V=\diag \left(\frac{\beta\gamma\varepsilon_1}{\alpha\psi_1},\ldots,
\frac{\beta\gamma\varepsilon_r}{\alpha\psi_r}\right) \oplus 0_{(p-r) \times (m-r)}
\in \F^{p\times m}.
\]
It follows from Lemma \ref{lem.div1}, \eqref{eq.deltai} and \eqref{eq.nui} that for $i=1,\ldots, r$
\[
\frac{\beta\gamma\varepsilon_i}{\alpha\psi_i}=\frac{\wh\beta_i \wt\varepsilon_i}{\wh\alpha_i\wt\psi_i},\quad
\gcd(\wh\beta_i ,\wt\psi_i)\doteq 1\quad \text{and}\quad \gcd(\wh\alpha_i,\wt\varepsilon_i)\doteq 1.
\]
Hence,  from this and using Lemma \ref{lem.div1}, Lemma \ref{lem.Fried} and \eqref{eq.chii}, it also holds for $i=1,\ldots, r$ that
\[
\frac{\beta\gamma\varepsilon_i}{\alpha\psi_i}=\frac{\wt\beta_i \wt\varepsilon_i}{\wt\alpha_i\wt\psi_i},\quad
\gcd(\wt\beta_i ,\wt\psi_i)\doteq 1,\quad  \gcd(\wt\alpha_i,\wt\varepsilon_i)\doteq 1\quad \text{and}\quad
\gcd(\wt\alpha_i,\wt\beta_i)\doteq 1.
\]
Now,  $\gcd(\varepsilon_i,\psi_i)\doteq 1$ and Lemma \ref{lem.Fried} imply  $\gcd(\wt\varepsilon_i,\wt\psi_i)\doteq 1$.
Thus,  taking into account that $\gcd(\wt\alpha_i,\wt\beta_i)\doteq 1$,  from Lemma \ref{lem.gcdprod} we get
$\gcd (\wt\beta_i \wt\varepsilon_i, \wt\alpha_i\wt\psi_i)\doteq 1$. We claim that
$\wt\beta_i\wt\varepsilon_i\mid\wt\beta_{i+1}\wt\varepsilon_{i+1}$ and
$\wt\alpha_{i+1}\wt\psi_{i+1}\mid \wt\alpha_{i}\wt\psi_{i}$ for $i=1,\ldots, r-1$. The claim implies that the Smith-McMillan form of $\wt{G}$ is
\be\label{eq.SMGt}
S_{\wt{G}}\doteq\diag \left(\frac{\wt\beta_1\wt\varepsilon_1}{\wt\alpha_1\wt\psi_1},\ldots,
\frac{\wt\beta_r\wt\varepsilon_r}{\wt\alpha_r\wt\psi_r}\right) \oplus 0_{(p-r) \times (m-r)}
\in \F^{p\times m}.
\ee
To prove the claim, we next observe some divisibility properties that will also be useful further on in this proof. By definition,
$\varepsilon_i\mid\varepsilon_{i+1}$ and then from \eqref{eq.deltai},
\be\label{eq.divdelta}
\delta_i\mid\delta_{i+1},\quad i=1,\ldots, r-1.
\ee
Applying Lemma \ref{lem_divnumden} with $a_1=\varepsilon_i$, $a_2=\varepsilon_{i+1}$ and $b_1=b_2=\alpha$,
we obtain
\be\label{eq.divteps}
\wt\varepsilon_i\mid\wt\varepsilon_{i+1},\quad i=1,\ldots,r-1.
\ee
Similarly, $\psi_{i+1}\mid \psi_i$ by definition and from \eqref{eq.nui},
\be\label{eq.divnu}
\nu_{i+1}\mid\nu_i,\quad i=1,\ldots, r-1.
\ee
By Lemma \ref{lem_divnumden} with $a_1=a_2=\beta\gamma$,  $b_1=\psi_i$ and $b_2=\psi_{i+1}$ we get
\be\label{eq.divtpsi}
\wt\psi_{i+1}\mid\wt\psi_i,\quad i=1,\ldots,r-1.
\ee
Finally, equations \eqref{eq.deltai}, \eqref{eq.nui}, \eqref{eq.divdelta} and \eqref{eq.divnu} collectively imply that, for all $i=1,\dots,r-1$,
$\wh\alpha_{i+1}\mid\wh\alpha_i$ and $\wh\beta_i\mid\wh\beta_{i+1}$. Applying once more Lemma \ref{lem_divnumden}, this time
with $a_1=\wh\beta_i$, $a_2=\wh\beta_{i+1}$,  $b_1=\wh\alpha_i$ and $b_2=\wh\alpha_{i+1}$, and taking
\eqref{eq.chii} into account,
we obtain
\be\label{eq.divtalpha}
\wt\alpha_{i+1}\mid\wt\alpha_i,\quad i=1,\ldots, r-1,
\ee
and
\be\label{eq.divtbeta}
\wt\beta_i\mid\wt\beta_{i+1},\quad i=1,\ldots, r-1.
\ee
The divisibility properties \eqref{eq.divteps}, \eqref{eq.divtpsi}, \eqref{eq.divtalpha} and \eqref{eq.divtbeta}
imply $\wt\beta_i\wt\varepsilon_i\mid\wt\beta_{i+1}\wt\varepsilon_{i+1}$ and
$\wt\alpha_{i+1}\wt\psi_{i+1}\mid \wt\alpha_{i}\wt\psi_{i}$ for $i=1,\ldots, r-1$, as claimed. 

Recall now that, in view of Definition \ref{rl_cop}, by construction the pair $(\alpha A, \gamma C)$ is right coprime and the pair $(\alpha A, \beta B)$ is left coprime. Hence, we can apply Theorem \ref{thm.1stproof} to $Q=\begin{bsmallmatrix} \alpha A & \beta B\\ \gamma C &
\frac{\beta\gamma}{\alpha}D\end{bsmallmatrix} \in\FR^{(n+p)\times (n+m)}$ and
$\wt{G}=\frac{\beta\gamma}{\alpha}D- \gamma C(\alpha A)^{-1}\beta B$ to conclude that, if $g$ is the largest index
in $\{1,\ldots, r\}$ such that $\wt\alpha_g\wt\psi_g\notin{\FR^\times}$
then the Smith form of $\alpha A$ is
\[
S_{\alpha A}\doteq I_{n-g}\oplus\diag\left(\wt\alpha_{g}\wt\psi_{g},\ldots, \wt\alpha_1\wt\psi_1\right).
\]
Now, from \eqref{eq.deltai} and \eqref{eq.chii} we get
\be\label{eq.delchi}
\frac{\wt\alpha_i}{\alpha}=\frac{1}{\delta_i\chi_i}.
\ee
Therefore, there are unimodular matrices $U_1,V_1\in\FR^{n\times n}$ such that
\[
U_1AV_1=\frac{1}{\alpha} I_{n-g}\oplus\diag\left(\frac{\wt\psi_g}{\delta_g\chi_g},\ldots, \frac{\wt\psi_1}{\delta_1\chi_1}
\right):
\]
we claim that the latter matrix is in Smith-McMillan form, yielding item i) in the statement. To prove the claim, note first that
$\wt\psi_{i+1}\mid\wt\psi_i$ for $i=1,\ldots, g-1$ because of  \eqref{eq.divtpsi} and $g\leq r$. Next,
 by \eqref{eq.delchi} $\delta_i\chi_i=\frac{\alpha}{\wt\alpha_i}$  and so
 $\delta_i\chi_i\mid\delta_{i+1}\chi_{i+1}$ follows from  \eqref{eq.divtalpha}; obviously, $\delta_g \chi_g \mid \alpha$ by \eqref{eq.delchi}. Finally we must verify that
 $\gcd(\wt\psi_i,\delta_i\chi_i)\eqdot 1$. By definition $\gcd(\psi_i,\varepsilon_i)\eqdot 1$ and, by \eqref{eq.nui},
 $\wt\psi_i\mid\psi_i$. Thus, Lemma \ref{lem.Fried} yields $\gcd (\wt\psi_i,\varepsilon_i)\eqdot 1$. But $\delta_i\mid \varepsilon_i$ by \eqref{eq.deltai}, and hence, again by Lemma \ref{lem.Fried} $\gcd(\wt\psi_i,\delta_i)\eqdot 1$. Consequently, it follows from Lemma \ref{lem.div3} that $\gcd (\wt\psi_i,\delta_i\chi_i)\eqdot
 \gcd (\wt\psi_i,\chi_i)$. Now, we had previously observed that $\gcd(\wt\psi_i,\wh\beta_i)\eqdot 1$, and since
$\chi_i\mid\wh\beta_i$ by \eqref{eq.chii} we conclude using Lemma \ref{lem.Fried} that $\gcd(\wt\psi_i,\chi_i)\eqdot 1$. Therefore, $\gcd (\wt\psi_i,\delta_i\chi_i)\eqdot \gcd (\wt\psi_i,\chi_i)\eqdot 1$,
 as claimed.

 Since the matrix $S_{\wt{G}}$ of \eqref{eq.SMGt} is a Smith-McMillan form of $\wt{G}$, we can use
 Theorem \ref{thm.1stproof} to obtain a Smith form of $Q = \beta \gamma P$:
 \[
 S_Q\doteq I_n\oplus\diag\left(\wt\beta_1\wt\varepsilon_1,\ldots, \wt\beta_r\wt\varepsilon_r\right)\oplus 0_{(p-r)\times(m-r)}.
 \]
 But $P=\frac{1}{\beta\gamma}Q$ and so $P$ is equivalent to
\[
 P_1=\frac{1}{\beta\gamma}I_n\oplus\diag\left(\frac{\wt\beta_1\wt\varepsilon_1}{\beta\gamma},\ldots,
 \frac{\wt\beta_r\wt\varepsilon_r}{\beta\gamma}\right)\oplus 0_{(p-r)\times(m-r)},
 \]
 Using \eqref{eq.nui} and \eqref{eq.chii},  we get $\frac{\wt\beta_i}{\beta\gamma}=\frac{1}{\chi_i\nu_i}$. Thus,
 \be\label{eq.SMP1}
 P_1=S_P=\frac{1}{\beta\gamma}I_n\oplus\diag\left(\frac{\wt\varepsilon_1}{\chi_1\nu_1},\ldots,
 \frac{\wt\varepsilon_r}{\chi_r\nu_r}\right)\oplus 0_{(p-r)\times(m-r)}
  \ee
is precisely the matrix appearing in item ii) of the statement.  It remains to prove that $S_P$ is in Smith-McMillan form. The argument is similar to the analogous task for item i), so we only sketch the procedure below without always explicitly declaring, in each step, which result in Appendix \ref{append} is needed.

By  \eqref{eq.divteps}, $\wt\varepsilon_i\mid\wt\varepsilon_{i+1}$, $i=1,\ldots, r-1$. Also, by
 \eqref{eq.divtbeta}, $\wt\beta_i\mid\wt\beta_{i+1}$, and so $\chi_{i+1}\nu_{i+1}=\frac{\beta\gamma}{\wt\beta_{i+1}}
 \mid \frac{\beta\gamma}{\wt\beta_i}=\chi_i\nu_i$, $i=1,\ldots, r-1$.  Since $\chi_1 \nu_1 \mid \beta \gamma$, it only remains to prove that $\gcd(\wt\varepsilon_i,
 \chi_i\nu_i)\eqdot 1$ for $i=1,\ldots, r$. In fact,  $\gcd (\wt\varepsilon_i,\psi_i)\eqdot 1$ follows from \eqref{eq.deltai}
and $\gcd(\varepsilon_i,\psi_i)\eqdot 1$. But, by \eqref{eq.nui}, $\nu_i\mid\psi_i$ and so
$\gcd(\wt\varepsilon_i,\nu_i)\eqdot 1$. By Lemma \ref{lem.div3}, $\gcd(\wt\varepsilon_i, \chi_i\nu_i)=
\gcd(\wt\varepsilon_i,\chi_i)$. Next, by \eqref{eq.deltai} $\wt\varepsilon_i=\frac{\varepsilon_i}{\delta_i}$ and $
\wh\alpha_i=\frac{\alpha}{\delta_i}$  and so $\gcd(\wt\varepsilon_i,\wh\alpha_i)\eqdot 1$. This and the
fact that, by \eqref{eq.chii}, $\chi_i\mid\wh\alpha_i$ imply $\gcd(\wt\varepsilon_i,\chi_i)\eqdot 1$, $i=1,\ldots, r$.

In conclusion, for $i=1,\ldots, r-1$, $\wt\varepsilon_i\mid\wt\varepsilon_{i+1}$, $\chi_{i+1}\nu_{i+1}\mid\chi_i\nu_i$,  $\chi_1 \nu_1 \mid \beta \gamma$,
and $\gcd(\wt\varepsilon_i, \chi_i\nu_i)\eqdot 1$ for $i=1,\ldots,r$. This concludes the proof. \qed \smallskip

\bigskip
Unlike Theorem \ref{thm.1stproof}, where the Smith-McMillan form of $G=D-CA^{-1}B$ determines the Smith
forms of $A$ and $P=\begin{bsmallmatrix} A & B\\C & D\end{bsmallmatrix}$ \textit{and vice versa}, Theorem
\ref{thm.1stproof-1} does not provide a straighforward way to obtain the Smith-McMillan form of $G$ out of
those of $A$ and $P=\frac{1}{\beta\gamma}\begin{bsmallmatrix} \alpha A & \beta B\\
\gamma C & \frac{\beta\gamma}{\alpha} D\end{bsmallmatrix}$.  This is our next goal. We need an auxiliary
result.

\begin{lem}\label{lem.SMFP}
Under the hypotheses of Theorem \ref{thm.1stproof-1}, let $P$ be the matrix as in \eqref{eq.firstdefP-1} and let
$U=\begin{bsmallmatrix}Y_{11}& Y_{12}\\Y_{21} & Y_{22}\end{bsmallmatrix}\in\FR^{(n+m)\times (n+m)}$ be
a unimodular matrix such that $\begin{bsmallmatrix} \alpha A & \beta B\end{bsmallmatrix}U=
\begin{bsmallmatrix}  I_n & 0\end{bsmallmatrix}$.  Then a Smith-McMillan form of $P$ is
\be\label{eq.SMFP2}
S_P\eqdot \frac{1}{\beta\gamma}I_n\oplus\diag\left(\frac{\nu_1}{\delta_1},\ldots,\frac{\nu_r}{\delta_r}\right)\oplus
0_{(p-r)\times(m-r)},
\ee
where $S_R\eqdot \diag\left(\frac{\nu_1}{\delta_1},\ldots,\frac{\nu_r}{\delta_r}\right)\oplus 0_{(p-r)\times(m-r)}$ is
a Smith-McMillan form of $R=\frac{1}{\beta}CY_{12}+\frac{1}{\alpha}DY_{22}$ and $n+r=\rank P$.
\end{lem}
\textbf{Proof}: Let $Q=\beta\gamma P \in \FR^{(n+p) \times (n+m)}$. As in the proof of Theorem \ref{thm.1stproof}, in particular \eqref{eq.1-firstprrosen},
$Q$ is equivalent to the matrix $I_n\oplus\left((\gamma C)Y_{12}+\left(\frac{\beta\gamma}{\alpha}D\right)Y_{22}\right)$.
Then $P=\frac{1}{\beta\gamma}Q$ is equivalent to $\frac{1}{\beta\gamma}I_n\oplus R$ and $\rank R=r$. Now,
$\beta\gamma R=(\gamma C)Y_{12}+ \frac{\beta\gamma}{\alpha}DY_{22}\in\FR^{p\times m}$ implying that
$\beta\gamma S_R\in\FR^{p\times m}$. This means that $\delta_1\mid\beta\gamma$ and so
$S_P\eqdot \frac{1}{\beta\gamma}I_n\oplus S_R$ as desired.\qed \smallskip

\begin{rem} \label{rem.detden} {\rm
A consequence of \cite[Theorem 2]{Coppel74} (see also the beginning of the proof of  \cite[Theorem 6]{Coppel74}) is that
the biggest denominator (in the sense of divisibility) of the Smith-McMillan form of any matrix over $\F$
is a least common denominator of the entries in that matrix. These theorems are proved in \cite{Coppel74}
assuming that $\F$ is the field of fractions of a PID $\FR$  but their proofs still hold when $\FR$ is a GCDD.
Hence, it follows from the above Lemma \ref {lem.SMFP}
that a least common denominator of $P$ is $\beta\gamma$.
}\end{rem}

\begin{thm} {\rm (Second part of Rosenbrock's theorem over fields of fractions of EDDs)}
\label{thm.1stproof-2}
Under the hypotheses of Theorem \ref{thm.1stproof-1}, let $P$ be the matrix as in \eqref{eq.firstdefP-1} and let the
matrix $S_P$ of \eqref{eq.SMFP2} be its Smith-McMillan form.  If  the Smith-McMillan form of $A$ is
\be\label{eq.SMFA2}
S_A\eqdot \diag\left(\frac{\beta_1}{\alpha_1},\ldots,\frac{\beta_n}{\alpha_n}\right),
\ee
then the Smith-McMillan form of $G=D-CA^{-1}B$ is
\be\label{eq.SMFG2}
S_G\eqdot\diag\left(\frac{\wt\alpha_n\wt\nu_1}{\wt\beta_n\wt\delta_1},\ldots,
\frac{\wt\alpha_{n-r+1}\wt\nu_r}{\wt\beta_{n-r+1}\wt\delta_r}\right)\oplus 0_{(p-r)\times(m-r)},
\ee
where $r=\rank P -n$,
\be\label{eq.mui}
\mu_i=\gcd(\alpha_{n-i+1},\delta_i),\quad \wt\alpha_{n-i+1}=\frac{\alpha_{n-i+1}}{\mu_i},\quad \wt\delta_i=\frac{\delta_i}{\mu_i},
\quad i=1, \ldots, r,
\ee
and
\be\label{eq.sigmai}
\sigma_i=\gcd(\beta_{n-i+1},\nu_i),\quad \wt\beta_{n-i+1}=\frac{\beta_{n-i+1}}{\sigma_i},\quad \wt\nu_i=\frac{\nu_i}{\sigma_i},
\quad i=1, \ldots, r.
\ee
\end{thm}
\textbf{Proof}:  Let $Q=\beta\gamma P \in \FR^{(n+p) \times (n+m)}$. Then, since $S_P$ as in \eqref{eq.SMFP2} is the Smith-McMillan form of $P$,
the Smith form of $Q$ is
\[
S_Q\eqdot I_n\oplus(\diag\left(\wh\delta_1 \nu_1, \ldots, \wh\delta_r\nu_r\right)\oplus 0_{(p-r)\times (n-r)},
\]
where $\wh\delta_i=\frac{\beta\gamma}{\delta_i}$, $i=1, \ldots, r$. Note that $\delta_{i+1}\mid\delta_i$ implies
$\wh\delta_i\mid\wh\delta_{i+1}$, $i=1, \ldots, r$. Similarly, since the Smith-McMillan form of $A$ is $S_A$ as in
\eqref{eq.SMFA2}, the Smith form of $\alpha A$ is
\[
S_{\alpha A}\doteq\diag\left( \wh\alpha_1\beta_1,\ldots, \wh\alpha_n\beta_n\right),
\]
where $\wh\alpha_i=\frac{\alpha}{\alpha_i}$, $i=1, \ldots, r$. Since $\alpha A\in\FR^{n\times n}$ and $\gcd(\alpha_1,\beta_1)
\eqdot 1$, $\alpha_1\mid\alpha$  by Lemma \ref{lem.div2}. Actually we know from Remark \ref{rem.detden}
that a least common denominator of $A$ is $\alpha_1$. Hence, $\alpha\eqdot \alpha_1$. Moreover, since $\alpha_i\mid\alpha_1$,
$\alpha_i\mid\alpha$ and $\wh\alpha_i\mid\wh\alpha_{i+1}$, $i=1,\ldots, r$, because $\alpha_{i+1}\mid\alpha_i$.

We can apply now the second part of Theorem \ref{thm.1stproof} to $Q = \begin{bsmallmatrix} \alpha A & \beta B\\ \gamma C &
\frac{\beta\gamma}{\alpha}D\end{bsmallmatrix}$ and $\alpha A$: If
$\wt{G}=\frac{\beta\gamma}{\alpha}D- \gamma C(\alpha A)^{-1}\beta B$
then its Smith-McMillan form is
 \[
 S_{\wt{G}}\eqdot\diag\left(\frac{\wh\delta_1 \nu_1}{\wh\alpha_n\beta_n},\ldots,
 \frac{\wh\delta_r \nu_r}{\wh\alpha_{n-r+1}\beta_{n-r+1}}\right)\oplus 0_{(p-r)\times(m-r)},
 \]
 where $\wh\alpha_1\beta_1\eqdot \cdots \eqdot \wh\alpha_{n-r}\beta_{n-r}\eqdot 1$ if $n\geq r$ and
 $\wh\alpha_i\beta_i\eqdot 1$ if $i<1$.
 Taking into account that $\wt{G}=\frac{\beta\gamma}{\alpha}G$, we can conclude that $G$ is equivalent to
 \[
 G_1=\diag\left(\frac{\alpha\wh\delta_1 \nu_1}{\beta\gamma\wh\alpha_n\beta_n},\ldots,
 \frac{\alpha\wh\delta_r \nu_r}{\beta\gamma\wh\alpha_{n-r+1}\beta_{n-r+1}}\right)\oplus 0_{(p-r)\times(m-r)}.
 \]
 By definition $\frac{\alpha}{\wh\alpha_i}=\alpha_i$ and $\frac{\beta\gamma}{\wh\delta_i}=\delta_i$. Thus
 \[
 G_1=\diag\left(\frac{\alpha_n\nu_1}{\delta_1\beta_n},\ldots,
 \frac{\alpha_{n-r+1} \nu_r}{\delta_r\beta_{n-r+1}}\right)\oplus 0_{(p-r)\times(m-r)}.
 \]
Now, $\alpha_n\nu_1\mid\cdots\mid \alpha_{n-r+1}\nu_r$ follows from $\alpha_n\mid\cdots \mid\alpha_1$ and $\nu_1\mid\cdots \mid\nu_r$, because these are the denominators and numerators  (different from 1, possibly)
of the diagonal entries of  $S_A$ and $S_P$, respectively.  Similarly, $\delta_r\beta_{n-r+1}\mid\cdots\mid \delta_1\beta_n$
because $\beta_1\mid\cdots \mid \beta_n$ and $\delta_r\mid\cdots \mid \delta_1$ are the numerators and denominators
(different from $\beta\gamma$, possibly) of the diagonal entries of $S_A$ and $S_P$, respectively.  Next,
$\gcd(\alpha_i,\beta_i)\eqdot 1$ for $i=1,\ldots, n$ and $\gcd(\nu_i,\delta_i)\eqdot 1$ for $i=1,\ldots, r$ because
$S_A$ of \eqref{eq.SMFA2} and $S_P$ of \eqref{eq.SMFP2} are the Smith-McMillan forms of $A$ and $P$, respectively.
Hence, applying Lemma \ref{lem.gcdprod} we get $\gcd(\alpha_{n-i+1}\nu_i, \delta_i\beta_{n-i+1})\eqdot
\gcd(\alpha_{n-i+1},\delta_i)\gcd(\nu_i,\beta_{n-i+1})$. But $\gcd(\alpha_{n-i+1},\delta_i)=\mu_i$ and
$\gcd(\nu_i,\beta_{n-i+1})=\sigma_i$ by definition (see \eqref{eq.mui} and \eqref{eq.sigmai}).
Also by \eqref{eq.mui} and \eqref{eq.sigmai},
$\wt\alpha_{n-i+1}=\frac{\alpha_{n-i+1}}{\mu_i}$, $\wt\delta_i=\frac{\delta_i}{\mu_i}$,
$\wt\beta_{n-i+1}=\frac{\beta_{n-i+1}}{\sigma_i}$ and $\wt\nu_i=\frac{\nu_i}{\sigma_i}$, $i=1, \ldots, r$.
Thus, for $i=1,\ldots,r$,
\[
\wt\alpha_{n-i+1}\wt\nu_i=\frac{\alpha_{n-i+1}\nu_i}{\mu_i\sigma_i}=
\frac{\alpha_{n-i+1}\nu_i}{\gcd(\alpha_{n-i+1}\nu_i,\beta_{n-i+1}\delta_i)},
\]
and
\[
\wt\beta_{n-i+1}\wt\delta_i=\frac{\beta_{n-i+1}\delta_i}{\mu_i\sigma_i}=
\frac{\beta_{n-i+1}\delta_i}{\gcd(\alpha_{n-i+1}\nu_i,\beta_{n-i+1}\delta_i)}.
\]
We have already seen that $\alpha_{n-i+1}\nu_i\mid\alpha_{n-i}\nu_{i+1}$ and $\delta_{i+1}\beta_{n-i}\mid
\delta_i\beta_{n-i+1}$. So we can apply Lemma \ref{lem_divnumden} with $a_i=\alpha_{n-i+1}\nu_i$
and $b_i=\delta_i\beta_{n-i+1}$ to get $\wt\alpha_{n-i+1}\wt\nu_i\mid \wt\alpha_{n-i}\wt\nu_{i+1}$
and $\wt\delta_{i+1}\wt\beta_{n-i}\mid\wt\delta_i\wt\beta_{n-i+1}$, $i=1,\ldots, r-1$. Therefore the matrix $S_G$
of \eqref{eq.SMFG2} is a Smith-McMillan form of $G$ and the theorem follows.\qed\smallskip

\begin{rem}{\rm
When $\FR=\K_{pr}(z)$ is the ring of proper rational functions and $A(z)\in\K[z]^{n\times n}$,
$B(z)\in\K[z]^{n\times m}$, $C(z)\in\K[z]^{p\times n}$ and $D(z)\in\K[z]^{p\times m}$, it is easy to check whether the condition
$\frac{\beta\gamma}{\alpha}D\in\FR^{p\times m}$ holds.
Recall that in this case $\alpha$, $\beta$ and $\gamma$ are $\frac{1}{z^{d_A}}$, $\frac{1}{z^{d_B}}$ and $\frac{1}{z^{d_C}}$
where $d_A$, $d_B$ and $d_C$ are, respectively, the degrees of the matrix polynomials $A(z)$, $B(z)$ and $C(z)$.
Then $\frac{\beta\gamma}{\alpha}D(z)\in\K_{pr}(z)^{p\times m}$ if and only if $d_B+d_C\geq d_A+d_D$ where
$d_D$ is the degree of $D(z)$.
Note that $P(z)=z^{d_B+d_C}\begin{bsmallmatrix}\frac{1}{z^{d_A}}A(z) & \frac{1}{z^{d_B}}B(z)\\
\frac{1}{z^{d_C}}C(z) &\frac{1}{z^{d_B+d_C-d_A}}D(z)\end{bsmallmatrix}$ and so the condition
$d_B+d_C\geq d_A+d_D$  guarantees that $Q(z) = \frac{1}{z^{d_B + d_C}} P(z)$ is a matrix of proper rational functions. Observe that
when $D(z) =0$ this degree inequality is always satisfied since $d_D = -\infty$. On the other hand,
\[
P\left(\frac{1}{z}\right)=\frac{1}{z^{d_B+d_C}}\begin{bmatrix} \rev A\left(z\right)&\rev B\left(z\right)\\\\
\rev C\left(z\right) &z^{d_B+d_C-d_A-d_D}\rev D\left(z\right)\end{bmatrix},
\]
where the reversal polynomials are defined in the standard way as $\rev A(z)= z^{d_A} A\left(\frac{1}{z}\right)$. Hence, taking into account that the
pole and zero structure at infinity of $G(z)\in\K(z)^{p\times m}$ is the pole and zero structure of $G\left(\frac{1}{z}\right)$
at $0$, Theorem \ref{thm.1stproof-2} says that, under the condition $d_B+d_C\geq d_A+d_D$, we can compute the pole and zero structure of $G\left(\frac{1}{z}\right)$
at $0$ by ``dividing'' the pole and zero structures of $P\left(\frac{1}{z}\right)$ and
$A\left(\frac{1}{z} \right) = \frac{1}{z^{d_A}}\rev A(z)$ at $0$ and ``simplifying numerators and denominators''.  Observe that the elementary divisors at infinity of $A(z)$ are the
elementary divisors at $0$ of $\rev A(z)$. If these are $ z^{e_1}$, \ldots , $z^{e_n}$ then, according to
\cite[Proposition 6.14]{AmMaZa14}, the invariant orders at infinity of $A(z)$ are $e_i-d_A$, $i=1,\ldots, n$. This property is revealed by the expression $A\left(\frac{1}{z}\right)=\frac{1}{z^{d_A}}\rev A(z)$.
}
\end{rem}

\begin{rem} \label{rem.G=A-1B} {\rm
When $G$ is written in the form $G=A^{-1}B$ ; i.e., $C=-I_n$ and $D=0_{n\times m}$, then $P$ becomes
$P=\frac{1}{\beta}\begin{bsmallmatrix} \alpha A & \beta B\\ I_n & 0\end{bsmallmatrix}$. We can easily annihilate
$\alpha A$ with a unimodular matrix:
\[
\begin{bmatrix} I_n & -\alpha A\\0& I_n\end{bmatrix}P=\begin{bmatrix} 0 & B\\\frac{1}{\beta}I_n & 0\end{bmatrix}.
\]
This means that $S_P \eqdot \frac{1}{\beta}I_n\oplus S_B$ where $S_B$ is a Smith-McMillan form of $B$ (recall that $\beta$
is a least common denominator of $B$).
It follows from Theorem \ref{thm.1stproof-1} that, under the assumptions made in that theorem and using the same notation,
\[
S_B\doteq\diag \left(\frac{\wt\varepsilon_1}{\nu_1\chi_1},\ldots,
\frac{\wt\varepsilon_r}{\nu_r\chi_r} \right) \oplus 0_{(n-r) \times (m-r)} \in \F^{n\times m}.
\]
So, in this case, Theorem \ref{thm.1stproof-1} is, as expected, a generalization of Theorem \ref{thm.1stpraux2} for the case when
$A\in\F^{n\times n}$ and $B\in\F^{n\times m}$, $G=A^{-1}B$ and $A$ and $B$ are left coprime matrices.

Conversely, if $A\in\F^{n\times n}$ and $B\in\F^{n\times m}$ are left coprime and
$S_B\eqdot\diag\left(\frac{\nu_1}{\delta_1},\ldots,\frac{\nu_r}{\delta_r}\right)\oplus 0_{(n-r)\times(m-r)}$
is a Smith-McMillan form of $B$ and $S_A$ of \eqref{eq.SMFA2} is that of $A$ then
$S_G$ of \eqref{eq.SMFG2} is a Smith-McMillan form of $G=A^{-1}B$. Let us note incidentally
that  the matrix $S_G$  of  \eqref{eq.SMFG2} is obtained from $S_A^{-1}S_B$ after reordering the
diagonal entries of $S_A$ and then simplifying as much as possible to get coprime numerators and denominators in each diagonal entry.

A similar result is obtained when $A\in\F^{n\times n}$ and $C\in\F^{p\times n}$ are right coprime and $G=CA^{-1}$.
}\end{rem}

\section{Conclusions} \label{sec.conclusions}We have extended Rosenbrock's theorem, proving it over general elementary divisor domains including those that are not principal ideal domains. Moreover, we have carefully studied the information delivered (or not) by a system matrix when the fundamental coprimeness conditions of Rosenbrock's theorem do not hold. For these analyses, we have generalized several results, known for rational and polynomial matrices, to general elementary divisor domains, and we have established some new auxiliary results. Our proof of Rosenbrock's theorem in this manuscript is not only valid under weaker assumptions on the base ring, but also more direct than other proofs available in the literature. Finally, we have explored how Rosenbrock's theorem can be extended to deal with matrices expressed as $G = D - C A^{-1} B$ when the matrices $A, B, C$ and $D$ have entries in the field of fractions of an elementary divisor domain. We therefore believe that our work can be of interest both to pure algebraists interested in general properties of elementary divisor domains and to practitioners who focus more concretely on rational and polynomial matrices, or on analytic and meromorphic matrices, or on integer matrices and matrices over the rational numbers.

\appendix

\section{Some divisibility properties over Greatest Common Divisor Domains}\label{append}

The notions of Greatest Common Divisor Domain and Elementary Divisor Domain are
defined, and their properties are studied, in many books,  e.g., the recent monographs \cite{Fried16, Shchedrykbook}.
Let us start by recalling the definitions of  Greatest Common Divisor Domain (GCDD) and B\'{e}zout Domains (BD):

\begin{description}

\item[\rm (a)] A GCDD is an integral domain $\FR$ such that any two elements $a,b\in\FR$  have a greatest common divisor (gcd) $d$;  i.e., $d\in \FR$ such that (1) $d$ is a divisor of both $a$ and $b$, or in formulae $d\mid a$ and $d\mid b$, (2) if $d'\in\FR$ satisfies $d'\mid a$ and $d'\mid b$, then $d'\mid d$. A more abstract characterization  is that, for every pair of elements $a,b$, there is a unique minimal principal ideal containing the ideal generated by $a$ and $b$. Moreover, any two elements of $\FR$ have a least common multiple (lcm). Over a GCDD, $\gcd(a,b)$ and
$\lcm(a,b)$ are related by the formula $ab\doteq \gcd(a,b)\lcm(a,b)$.

\smallskip

\item[\rm (b)] A BD is an integral domain $\FR$  such that  any two elements $a,b\in\FR$ have a greatest common divisor, $\gcd(a,b)$,  satisfying $pa+qb \doteq \gcd(a,b)$ for some $p,q\in\FR$.  Equivalent characterizations are that every finitely generated ideal is principal, or that the sum of any two principal ideals is principal.

\end{description}

Thus, a BD is always a GCDD. More generally, the class inclusions of Figure \ref{fig:1} hold  {
\cite{Fried16, Shchedrykbook} (see \cite[p. 226]{hel43} for the strict inclusion $\text{PID}\subset\text{AD}$ and \cite[pp. 159, 161-162]{henrik55} for $\text{AD}\subset\text{EDD}$)}. To our knowledge and at the time of writing the present manuscript, it is still an open question to determine whether the class inclusion EDD $\subseteq$ BD is strict  or not; see \cite{BovDi19}, for example.
\begin{figure}[h]
\begin{tikzpicture}
\node(o) at (12.5,0){};
\node (i) at (14,0) {ID};
\node(gcd) at (15.5,0) {GCDD};
\node (uf) at (18.5,0.85) [draw,minimum width=4.5cm,minimum height=0.5cm]{UFD};
\node(b) at (17,0)  {BD};
\node(ed) at (18.55, 0) {EDD};
\node(ad) at (20, 0) {AD};
\node(pi) at (21.5, 0) {PID};
\node(e) at  (23,0) {ED};
\node(f) at (24.5, 0) {Fields};
\node at($(i)!0.45!(gcd)$) {$\supset$};
\node at (16,0.45){\rotatebox{45}{$\supset$}};
\node at ($(gcd)!0.6!(b)$) {$\supset$}; 
\node at ($(b)!0.45!(ed)$) {$\supseteq$};
\node at ($(ed) !0.45!(ad)$){$\supset$};
\node at (21,0.45) {\rotatebox{315}{$\supset$}};
\node at ($(ad)!0.45!(pi)$) {$\supset$};
\node at ($(pi)!0.45!(e)$) {$\supset$};
\node at ($(e)!0.45!(f)$) {$\supset$};
\end{tikzpicture}
\[
\begin{array}{l}
\text{ID}=\text{Integral Domains }, \; \text{GCDD}=\text{ Greatest Common Divisor Domains},\\
\text{BD}= \text{B\'{e}zout Domains},\;\text{EDD}= \text{Elementary Divisor Domains},\\
{\text{AD}=\text{Adequate Domains},}\; \text{UFD}= \text{Unique Factorization Domains},\\
\text{PID}=\text{ Principal Ideal Domains},\; \text{ED}= \text{Euclidean Domains}
\end{array}
\]
\caption{Inclusion relations between different classes of  commutative rings discussed in this paper. All inclusions in the figure are known to be strict with the exception of the one between B\'{e}zout domains and elementary divisor domains}
\label{fig:1}
\end{figure}

Several important classes of rings have not been depicted in Figure \ref{fig:1}; in particular non-commutative rings. Although non-commutative adequate and elementary divisor rings have been
studied (see \cite{BovDi22,GatZaba99} and the references therein), they are beyond the scope of this article. We leave to possible future research the question of whether the ideas and results of this manuscript still hold for matrices with entries in those rings.

In the rest of this appendix we give an account of the divisibility properties that are used in the manuscript. Although in the manuscript we focus on EDDs, these properties only require the assumption that $\FR$ is a GCDD.
In the literature, these properties are often stated under the assumption that the underlying ring is a BD. However, their proofs for elements in any GCDD seem to be more difficult to find. For readers' convenience, we present them as a series of lemmas with corresponding proofs.
We will provide references for those results that we know can be found
elsewhere.

\begin{lem}\cite[Theorem 49]{kaplansky}\label{lem.div1} Let $\FR$ be a GCDD.
If $a,b,c\in\FR$ then $\gcd(ac,bc)\doteq c\gcd(a,b)$.
\end{lem}
\textbf{Proof}: Let $d_1\doteq \gcd(a,b)$ and $d_2\doteq \gcd(ac,bc)$. Then $d_1\mid a$, $d_1\mid b$ and so
$cd_1\mid\gcd(ac,bc)$ implying  $cd_1\mid d_2$. It follows then that there is a multiple of $d_1$, say, $g$, such that
$d_2=cg$. Since $d_2\mid ac$ and $d_2\mid bc$, and taking into account that $\FR$ is an integral domain, it must be $g\mid a$ and $g\mid b$.
But then $g\mid\gcd(a,b)\doteq d_1$. As $g\mid d_1\mid g$ and $d_2=cg$ we conclude $d_2\doteq cd_1$.\qed\smallskip

\begin{lem} \label{lem.div1lcm}  Let $\FR$ be a GCDD.
If $a,b,c\in\FR$ then $\lcm(ac,bc)\doteq c\lcm(a,b)$.
\end{lem}
\textbf{Proof}: The statement is an immediate corollary of Lemma \ref{lem.div1} and the formula $ab\doteq\gcd(a,b)\lcm(a,b)$.\qed\smallskip

\begin{lem}\label{lem.div2}
Let $\FR$ be a GCDD. If $a,b,c\in\FR$ satisfy $\gcd(a,b)\doteq 1$ and $a\mid bc$ then $a\mid c$.
\end{lem}
\textbf{Proof}:  By assumption $a\mid bc$ and it is plain that $a\mid ac$. Then
$a\mid\gcd(ac,bc)$ and, by Lemma \ref{lem.div1}, $\gcd(ac,bc)\doteq c\gcd(a,b)\doteq c$.\qed\smallskip

\begin{lem}\label{lem_divnumden}
Let $\FR$ be a GCDD and let $a_1,a_2,b_1,b_2$ be elements of $\FR$ such that $a_1\mid a_2$ and $b_2\mid b_1$.
Then
\be\label{eq.divgcd}
\frac{a_1}{\gcd(a_1,b_1)}\;\Big|\;\frac{a_2}{\gcd(a_2,b_2)} \quad \text{ and } \quad
\frac{b_2}{\gcd(a_2,b_2)}\;\Big|\;\frac{b_1}{\gcd(a_1,b_1)}.
\ee
\end{lem}
\textbf{Proof}:
Let $d_1\doteq\gcd(a_1,b_1)$ and $d_2\doteq\gcd (a_2,b_2)$. From this and  $a_1\mid a_2$, $b_2\mid b_1$ it follows that
there are elements $x,y,p_1,p_2,q_1,q_2\in\FR$ such that $b_1=xb_2$, $a_2=ya_1$, $a_1=p_1d_1$,
$a_2=p_2d_2$, $b_1=q_1d_1$ and $b_2=q_2d_2$. Then
\[
\begin{array}{rcl}
p_2q_1d_1&\stackrel{b_1=q_1d_1}{=}&p_2b_1\stackrel{b_1=xb_2}{=} p_2xb_2\stackrel{b_2=q_2d_2}{=}p_2xq_2d_2\\
&\stackrel{a_2=p_2d_2}{=}&xq_2 a_2\stackrel{a_2=ya_1}{=}xq_2ya_1\stackrel{a_1=p_1d_1}{=}xq_2yp_1d_1.
\end{array}
\]
Since $\FR$ is an integral domain $p_2q_1=xq_2yp_1$ and so $p_1\mid p_2q_1$. Now, $d_1\doteq\gcd(a_1,b_1)\doteq
\gcd(p_1d_1,q_1d_1)$ and by  Lemma \ref{lem.div1}, $d_1\doteq d_1\gcd(p_1,q_1)$. This means that
$\gcd(p_1,q_1)\doteq 1$. Since $p_1\mid p_2q_1$ and $\gcd(p_1,q_1)\doteq 1$, it follows from
Lemma  \ref{lem.div2} that $p_1\mid p_2$. That is, $\frac{a_1}{d_1}\;\big|\; \frac{a_2}{d_2}$.

Similarly, it follows from $p_2q_1=xq_2yp_1$ that $q_2\mid p_2q_1$, and from
$d_2\doteq \gcd(a_2,b_2) \doteq \gcd(p_2d_2,q_2d_2) \doteq  d_2 \gcd(p_2,q_2)$ that $\gcd(p_2,q_2) \doteq 1$. So
$q_2\mid q_1$ and, equivalently, $\frac{b_2}{d_2}\;\big|\; \frac{b_1}{d_1}$ and \eqref{eq.divgcd} follows.\qed \smallskip

\begin{lem}\label{lem.gcdprod}
Let $\FR$ be a GCDD and let $a,b,c,d$ be elements of $\FR$ such that $\gcd(a,b)\doteq 1$ and $\gcd(c,d)\doteq 1$.
Then $\gcd(ac,bd)\eqdot \gcd(a,d)\gcd(b,c)$.
\end{lem}
\textbf{Proof}:
To simplify notation, define $g \eqdot \gcd(ac,bd)$, $f \eqdot \gcd(a,d)$ and $h \eqdot \gcd(b,c)$.  It is clear that the product $fh$ divides both $ac$ and $bd$, and hence $fh \mid g$.

For the reverse divisibility property, define
$e_1\eqdot\gcd(a,g)$ and $e_2\eqdot \gcd(b,g)$. Observe (Lemma \ref{lem.div1}) that $\gcd\left(\frac{a}{e_1},\frac{g}{e_1}\right)\eqdot 1$ and $\gcd\left(\frac{b}{e_2},\frac{g}{e_2}\right)\eqdot 1$. It follows from $g\mid ac$ that $\frac{g}{e_1}\mid\frac{a}{e_1}c$. Since $\gcd\left(\frac{a}{e_1},\frac{g}{e_1}\right)\eqdot 1$,
we conclude (Lemma \ref{lem.div2}) that $\frac{g}{e_1}\mid c$. Similarly $\frac{g}{e_2}\mid d$.
But, by assumption,  $\gcd (c,d)\eqdot 1$. Thus  $\gcd\left(\frac{g}{e_1},\frac{g}{e_2}\right)\eqdot 1$ which implies (Lemma \ref{lem.div1}) that
$g\gcd(e_1,e_2)\eqdot\gcd(ge_2,ge_1)\eqdot e_1e_2$. On the other hand,
\[ \gcd(e_1,e_2)\eqdot\gcd(\gcd(a,g),\gcd(b,g))\eqdot \gcd(a,b,g)\eqdot \gcd(\gcd(a,b),g) \eqdot \gcd(1,g) \eqdot 1, \] implying $g\eqdot e_1e_2$. At this point, observe that $e_1 \mid a$ by definition and that by the argument above $e_1 \doteq \frac{g}{e_2} \mid d$. Hence, $e_1 \mid f$. We can show similarly that $e_2 \mid h$, and hence $g \doteq e_1 e_2 \mid fh$.
\qed \smallskip

Straightforward consequences of Lemma \ref{lem.gcdprod} are the following results.
\begin{lem}\label{lem.div3}
Let $\FR$ be a GCDD. If $a,b,c\in\FR$ and $\gcd(a,b)\doteq 1$ then $\gcd(a,bc)\doteq \gcd (a,c)$.
\end{lem}

\begin{lem}\label{lem.Fried}
Let $\FR$ be a GCDD. Let $a,b,c\in\FR$. Then $\gcd(a,b)\doteq 1$ and $\gcd (a,c)\doteq 1$
if and only if $\gcd(a,bc) \doteq 1$.
\end{lem}

For the ``only if'' part in Lemma \ref{lem.Fried} observe that if $\gcd(a,bc)\doteq 1$ then $\gcd(a,b) \mid \gcd(a,bc) \doteq 1$. So, $\gcd(a,b) \doteq 1$.
Similarly, $\gcd (a,c) \doteq 1$.

\section{Proof of Proposition \ref{teo_copr_an}} \label{appendcoprimeness}

\noindent i) $\Rightarrow$ ii) Let $\left[\begin{smallmatrix}
S \\
0
\end{smallmatrix} \right]$ be any arbitrary Smith form of $\left[\begin{smallmatrix}
  G_1 \\
  G_2
\end{smallmatrix} \right] $. Then, it holds
\begin{equation}\label{eq.1appendcopri}
\begin{bmatrix}
  G_1 \\
  G_2
\end{bmatrix} = \begin{bmatrix}
                  W_{11} & W_{12} \\
                  W_{21} & W_{22}
                \end{bmatrix}
                \begin{bmatrix}
                  S \\
                  0
                \end{bmatrix} Z,
\end{equation}
where
$\left[ \begin{smallmatrix}
                  W_{11} & W_{12} \\
                  W_{21} & W_{22}
                \end{smallmatrix} \right] \in \FR^{(p+q)\times (m+(p+q-m))}$ is unimodular,
$S\in \FR^{m\times m}$ is diagonal, and $Z \in \FR^{m\times m}$ is unimodular. Then
\[
G_1 = W_{11} (SZ) \quad \mbox{and} \quad G_2 = W_{21} (SZ)
\]
and hence $SZ$ is a common right divisor of $G_1$ and $G_2$. Therefore $SZ$ is unimodular. Hence $S$ is also unimodular and $S \doteq I_m$.

\bigskip
\noindent ii) $\Rightarrow$ iii) Let us write \eqref{eq.1appendcopri} with $S = I_m$ as
\[
\begin{bmatrix}
  G_1 \\
  G_2
\end{bmatrix} = \left(\begin{bmatrix}
                  W_{11} & W_{12} \\
                  W_{21} & W_{22}
                \end{bmatrix}
                \begin{bmatrix}
                  Z & 0 \\
                  0 & I
                \end{bmatrix} \right)
                \begin{bmatrix}
                  I_m \\
                  0
                \end{bmatrix} = U^{-1} \begin{bmatrix}
                  I_m \\
                  0
                \end{bmatrix} \; \mbox{with} \; U^{-1} = \begin{bmatrix}
                  W_{11} & W_{12} \\
                  W_{21} & W_{22}
                \end{bmatrix}
                \begin{bmatrix}
                  Z & 0 \\
                  0 & I
                \end{bmatrix} .
\]
Then $U^{-1}$ is unimodular and $U$ is the unimodular matrix in the statement of Proposition \ref{teo_copr_an}-iii).

\bigskip
\noindent iii) $\Rightarrow$ iv) The equation
$\begin{bsmallmatrix} G_1 \\ G_2 \end{bsmallmatrix} = U^{-1} \begin{bsmallmatrix} I_m\\ 0 \end{bsmallmatrix}$ and the fact that $U^{-1}$ is unimodular imply that
$U^{-1} = \begin{bsmallmatrix} G_1 & C \\ G_2 & D \end{bsmallmatrix}$ is unimodular for some matrices $C\in\FR^{p\times (p+q-m)}$, $D\in\FR^{q\times (p+q-m)}$.

\bigskip
\noindent iv) $\Rightarrow$ v) Consider the unimodular matrix
$
\begin{bsmallmatrix} X & Y \\ E & F \end{bsmallmatrix} = \begin{bsmallmatrix} G_1 & C \\ G_2 & D \end{bsmallmatrix}^{-1} \, ,
$
with $X\in\FR^{m\times p}$, $Y\in\FR^{m\times q}$. Then
$
\begin{bsmallmatrix} X & Y \\ E & F \end{bsmallmatrix}  \begin{bsmallmatrix} G_1 & C \\ G_2 & D \end{bsmallmatrix} = \begin{bsmallmatrix} I_m & 0 \\ 0 & I \end{bsmallmatrix}
$
implies $XG_1+YG_2=I_m$.

\bigskip
\noindent v) $\Rightarrow$ i) If $F \in \FR^{m\times m}$ is a common right divisor of $G_1$ and $G_2$ and $G_1 = \widetilde{G_1} F$ and $G_2 = \widetilde{G_2} F$, then $(X \widetilde{G_1} + Y \widetilde{G_2})F = I_m$. This implies $\det(X \widetilde{G_1} + Y \widetilde{G_2}) \, \det F = 1$. Therefore $\det F \in {\FR^\times}$ and $F$ is unimodular.

\end{document}